\def\**{$*\!*\!*$}


 at 12truept
\font\tenmsy=msym10
\textfont8=\tenmsy
\mathchardef\ssm="7872

\input psfig
\mathsurround = 2pt
\abovedisplayskip=6pt
\belowdisplayskip=6pt

\def \QED {\rlap{$\sqcup$}$\sqcap$\smallskip}

\def\ref{\hangindent=1pc \hangafter=1 \noindent}

\def\IMSmarkvadjust{0 pt}
\def\IMSmarkhadjust{0 pt}
\def\IMSmarkhpadding{0 pt}
\def\SBIMSMark#1#2#3{
 \font\SBF=cmss10 at 10 true pt
 \font\SBI=cmssi10 at 10 true pt
 \setbox0=\hbox{\SBF \hbox to \IMSmarkhpadding{\relax}
                Stony Brook IMS Preprint \##1}
 \setbox2=\hbox to \wd0{\hfil \SBI #2}
 \setbox4=\hbox to \wd0{\hfil \SBI #3}
 \setbox6=\hbox to \wd0{\hss
             \vbox{\hsize=\wd0 \parskip=0pt \baselineskip=10 true pt
                   \copy0 \break%
                   \copy2 \break%
                   \copy4 \break}}
 \dimen0=\ht6   \advance\dimen0 by \vsize \advance\dimen0 by 8 true pt
                \advance\dimen0 by -\pagetotal
	        \advance\dimen0 by \IMSmarkvadjust
 \dimen2=\hsize \advance\dimen2 by .25 true in
	        \advance\dimen2 by \IMSmarkhadjust

%
%
  \openin2=publishd.tex
  \ifeof2\setbox0=\hbox to 0pt{}
  \else 
     \setbox0=\hbox to 3.1 true in{
                \vbox to \ht6{\hsize=3 true in \parskip=0pt  \noindent  
                {\SBI Published in modified form as part of}\hfil\break
``Dynamics of quadratic polynomials, I-II'',
{\it  Acta Math.}~{\bf 178} (1997), 185-297.
 
                \vfill}}
  \fi
  \closein2
  \ht0=0pt \dp0=0pt
 \ht6=0pt \dp6=0pt
 \setbox8=\vbox to \dimen0{\vfill \hbox to \dimen2{\copy0 \hss \copy6}}
 \ht8=0pt \dp8=0pt \wd8=0pt
 \copy8
 \message{*** Stony Brook IMS Preprint #1, #2. #3 ***}
}

\SBIMSMark{1993/9}{October 1993}{}

\centerline{\bf Geometry of quadratic polynomials:
 moduli, rigidity and local connectivity.}
\smallskip

\medskip

\centerline{Mikhail Lyubich\footnote{*}{\rm Supported in part by 
NSF grant DMS-8920768 and a Sloan Research Fellowship.}}\smallskip
\centerline{Mathematics Department and IMS, SUNY Stony Brook}\bigskip

\centerline {\bf \S1. Introduction.}\bigskip

A key problem in holomorphic dynamics is to classify
 complex quadratics
$z\mapsto z^2+c$
 up to topological conjugacy. The Rigidity Conjecture would assert that any
non-hyperbolic polynomial  is {\sl topologically rigid},
 that is, not topologically conjugate
to any other polynomial. This would imply density of hyperbolic 
polynomials in the complex quadratic family (Compare Fatou [F, p. 73]).
A stronger conjecture  
usually abbreviated as MLC would assert that the Mandelbrot set is locally
connected (see [DH1]).
  
A while ago MLC was proven for quasi-hyperbolic points by Douady and Hubbard,
and for boundaries of hyperbolic components by Yoccoz.
More recently Yoccoz proved MLC for all at most finitely renormalizable
parameter values (see [H], [M2] for the exposition of this work and 
closely related work of Branner and Hubbard  [BH] on rigidity of cubics). One of
our goals is to prove MLC for some infinitely renormalizable parameter values.
Loosely speaking, we need all renormalizations to have bounded combinatorial 
rotation number (assumption C1) and sufficiently high
  combinatorial type (assumption C2)
(see \S 2 for the  precise statement of the assumptions).

This result is based on a complex version of a theorem of [L2] which
says that the scaling factors characterizing the geometry
of a real non-renormalizable quasi-quadratic map  decay exponentially. 
Its complex counterpart proved below  (Theorem I) says that  the moduli of the 
principle nest of annuli  grow  linearly (this result does not need
any a priori assumptions). This makes finitely renormalizable  maps
{\sl geometrically tame}
in the sense that the return maps are becoming purely quadratic in small
scales.  In the infinitely renormalizable case satisfying
 assumptions (C1) and (C2) Theorem I implies
 complex a priori bounds  (Theorem II)
  (that is, the  bounds from below for the moduli of the fundamental 
annuli of $R^nf$). 

For  real quadratic polynomials of bounded  combinatorial 
type the complex a priori bounds
were obtained by Sullivan [S]. Our 
result complements the Sullivan's result in the unbounded case.
Moreover, it gives a background for Sullivan's renormalization theory
for some bounded type polynomials outside the real line 
where the problem of a priori
bounds was not handled before for any single polynomial. 

An important consequence of a priori bounds is absence of invariant 
measurable line fields on the Julia set (McMullen [McM])
 which is equivalent to quasi-conformal (qc) rigidity. To prove 
stronger topological
rigidity we construct a qc conjugacy between any two
topologically conjugate polynomials (Theorem III). We do this by 
means of a pull-back argument, based on the linear growth of moduli
and a priori bounds. Actually the argument gives  the stronger
 {\sl combinatorial rigidity}  which implies MLC.

Local connectivity of
the Julia set  is also a general consequence of a priori bounds 
(see Hu and Jiang [HJ], [J]),
so we have it under assumptions (C1) and (C2). Note that 
  Douady and Hubbard gave an example of an infinitely
renormalizable polynomial with non-locally connected Julia set
 (see Milnor's version of the example in [M2]).
 In this example the combinatorial rotation numbers of
the fixed points of $R^m f$ are highly unbounded, which is ruled out by
our first assumption.

We complete the paper with an application to the real quadratic family.
Here we can give a precise dichotomy (Theorem IV): on each renormalization
level we either observe a big modulus, or essentially bounded geometry.
This allows us to combine the above considerations with Sullivan's
argument for bounded geometry case, and to obtain a new proof of the
rigidity conjecture on the real line (compare McMullen [McM] and Swiatek [Sw]).  

This paper is organized as follows.
\S 2 contains a combinatorial framework: Yoccoz puzzle, principle nest,
usual and generalized renormalization.
Theorems I and II on geometric moduli in the 
dynamical plane are proved in \S 3, 
Theorem III yielding MLC  under the above assumptions is proved in \S 4. 
 The real case is discussed in \S 5.

\medskip\noindent {\bf Acknowledgement.} I would like to thank Jeremy Kahn
for a fruitful suggestion (see \S 3), and Curt McMullen
and Mitsuhiro Shishikura for useful comments on the results. 
I also thank Scott Sutherland and Brian Yarrington for making the 
computer pictures.
This work was done in June 1993 during the Warwick Workshop on
hyperbolic geometry. I am grateful to 
the organizers, particularly  David Epstein and Caroline Series, 
for that wonderful time.

\bigskip

\centerline{\bf \S 2. Principle nest and renormalization.}\bigskip

We refer to [D1], [DH2] and [M1] for the background in polynomial-like 
mappings and tuning (which is called below
 ``quadratic-like renormalization"),
and to [H] and [M2] for the introduction to the Branner-Hubbard-Yoccoz puzzle.

Let  $f: U'\rightarrow U$ be a quadratic-like map  with
connected Julia set. Let us start with an appropriate combinatorially defined
$c$-symmetric domain (``puzzle-piece")
$V^0\subset U'$  such that 
(see the construction below). 
Let us then consider the first return of the critical point back to $V^0$,
and pull $V^0$ back along the corresponding piece of the critical orbit.
In such a way we obtain the critical puzzle-piece $V^1\subset V^0$ of the
first  level.
 If we do the same replacing $V^0$ by $V^1$, we 
obtain the critical puzzle-piece $V^2\subset V^1$ of the second level. 
Proceeding in this  manner we will construct the
{\sl principle nest}
$$V^0\supset V^1\supset V^2\ldots$$
of puzzle-pieces. It may happen that on some level the quadratic-like map
$f^{\pi(t)}: V^{t+1}\rightarrow V^t$  has a connected Julia set
(which is equivalent to having non-escaping critical point). 
Then we say that $f$ is {\sl q-renormalizable}, or
that $f$ admits the {\sl quadratic-like renormalization} 
$R f=f^{\pi(t)}: V^{t+1}\rightarrow V^t$
(usually such a map is just called ``renormalizable" but we need 
 to distinguish  the quadratic-like renormalization from
  the generalized renormalization defined below).
 In this case
the puzzle-pieces $V^n$ shrink down to the Julia set $J(R f)$.
Otherwise by the
Yoccoz Theorem they shrink down to the critical point.



Let now $V^n_i\subset V^{n-1}$ denote the pull-backs of $V^{n-1}$
corresponding to the first returns of the points  
$x\in\omega(c)\cap V^{n-1}$ back to $V^{n-1}$, 
numbered so that $V^n\equiv V^n_0$.
 The first return map
$$g_n : \cup V^n_i\rightarrow V^{n-1}$$
we call the $n$-fold (generalized) {\sl renormalization} of $f$. 
If  $f$ admits a quadratic-like renormalization,
 then the return time to $V^{n-1}$ is uniformly
bounded on $\omega (c)$, and hence the domain of the renormalized
maps consists of only finitely many components $V^n_i$.
We select the initial puzzle-piece $V^0$ in such a way that
$V^1$ is compactly contained in $V^0$ (see Lemma 0 below). 
Then the puzzle-pieces $V^n_i,\; n\geq 2,$ are
 compactly contained in $V^{n-1}$ as well, and hence 
 the $g_n$ are 
{\sl generalized polynomial-like maps}
in the sense of [L1].

 Let us call this return to level $n-1$ ``central" if $g_n c\in V^n$.
If we have several subsequent central returns, we refer to a
{\sl cascade} of central returns. In the q-renormalizable case the
sequence of puzzle-pieces $V^n$ ends with an infinite cascade of central
returns.  Let us denote by $\kappa=\kappa(f)$ the number of the levels
on which the  non-central return occurs. The map $f$ 
 admits a quadratic-like renormalization 
iff $\kappa<\infty$.  

Let us now construct the initial puzzle-piece $V^0$.  Let $q/p$ be the
combinatorial rotation number of the dividing fixed point $\alpha$.
This means that
there are $p$ disjoint {\sl external rays} $\Gamma_i$
 landing at $\alpha$ which are permuted by the 
dynamics with rotation number $q/p$.  They cut $U'$
 into $p$ initial {\sl Yoccoz puzzle-pieces}. 
(Warning: 
we define the external rays in a
non-canonical way via a conjugacy to a polynomial. These rays are not
necessarily the external rays for the original map $f$, whose
geometry we would not be able to control).

Let $\Omega^0$ be  the  critical puzzle-piece, that is, the one
containing the critical point. 
Let us pull it back along the critical orbit in the same way as we did above
with $V$-pieces.
Then in the beginning we may observe a central cascade 
$$\Omega^0\supset \Omega^1\supset \ldots$$
In what follows we always assume that this first cascade is finite,
that is, there is an $N$, such that
$f^p c\in \Omega^{N-1}\backslash \Omega^N$
(this can be viewed as a part of  Assumption (C2) that the combinatorial type
is sufficiently high). 


Let now $\gamma$ and $\gamma'=-\gamma$ be the periodic and co-periodic points of
period $p$ belonging to $\Omega^0$. Let us truncate $\Omega^N$ 
by the external rays
landing at these points. The critical puzzle-piece obtained in such
a way is the desired $V^0$ (see Figure 1).

\proclaim Lemma 0. The puzzle-piece $V^1$ is compactly contained in $V^0$.

\noindent{\bf Proof.} The argument below is not the shortest possible,
but it will later give us important extra information (see Lemma 2). 
Let $W_i,\; i=1,...,p-1,$
 be the puzzle-pieces bounded by the
external rays landing at $\alpha'$ and the equipotential $f^{-1} U$,
numbered in such a way that $f^i W_i\supset U_0$.

Take a point $z\in \Omega^0\cap J(f)$, push it forward by iterates of $g=f^p$,
 and find the first  moment $r=r(z)$ (if any) such that $g^rz$ lands
either at $W_i$ (where $i=i(z)$) or at $V^0$.
 In the first case consider the pull-back $X(z)$ of $W_i$ along the 
orbit of $z$, in the second case consider the pull-back $Y(z)$ of $V^0$.
(The points which are not covered by the sets $X_j$ and $Y_j$
 form an invariant Cantor set
in $\Omega^N\backslash V^0$.)

Let us now define a map $G: X_j\cup Y_j\rightarrow U^0\cup V^0$ in the following
way:  $G=f^i\circ g^r$ on $X_j$ and $G=g^r$ on $X_j$. Then every $X_j$ is
univalently mapped onto $\Omega^0$, while $Y_j$ is univalently mapped onto
$V^0$. Let us push $c$ forward by iterates of $G$ until the first moment $l$ it
lands at a set $Y_j$. Let $Q\ni c$ be the pull-back of $\Omega^0$ under $G^l$.
 Then  
  $$G^l : (Q, V^1)\rightarrow (\Omega^0, Y_j). \eqno (0)$$ 
As $Y_j$ is compactly contained in $\Omega^0$, $V^1$ is compactly contained 
in $Q$. Observe finally that $Q\subset V^0$, since $Q$ may not intersect the
boundary of $V^0$.  \QED

\smallskip
If $f$  infinitely q-renormalizable
 then we can repeat the above construction 
on the corresponding quadratic-like levels, and consider
 the full {\sl canonical nest} of puzzle-pieces:
$$ V^{0,0}\supset V^{0,1}\supset\ldots V^{0,t(0)+1}\supset
  V^{1,0}\supset V^{1,1}\supset\ldots\supset V^{1,t(1)+1}\supset V^{2,0}
  \ldots $$
Here the first index counts the quadratic-like  levels, while
the second one  counts the levels in between. The maps
$V^{m,t(m)+1}\rightarrow V^{m,t(m)}$ are quadratic-like with non-escaping
critical point, while $V^{m+1,0}$ are the critical puzzle-pieces obtained by
the above procedure applied to these maps. (The choice of the cutting level
$t(m)$ is not  canonical). We skip the first index
when we work in between two quadratic-like levels.

Let us finish this section with specifying exact conditions under which
we will prove MLC. To these end we need several notions. A {\sl limb} of
the Mandelbrot $M$ set is the connected components of $M\backslash \{c_0\}$
(which does not contain 0)
where $c_0$ is a bifurcation point on the main cardioid. A limb is specified by 
specifying a combinatorial rotation number at the dividing fixed point.
If we remove from a limb a neighborhood of its root $c_0$,
 what is left  we call a {\sl truncated} limb.
By a {\sl (truncated) secondary limb } we mean
the similar object corresponding to the second bifurcation
from the main cardioid ( see Figure~2). 

Two  quadratic-like maps are called {\sl hybrid (or internal) equivalent}
if they are conjugate by a qc map $h$ with $\bar\partial h=0$ almost
everywhere on the Julia set. By the Douady-Hubbard  Straightening Theorem
[DH2], any hybrid class with connected Julia set contains a unique
quadratic polynomial $z\mapsto z^2+c$. So such hybrid classes are labeled
by points on the Mandelbrot set. 

Let ${\cal F}$ denote the class of maps admitting 
infinitely many  quadratic-like renormalizations and
  satisfying the following assumptions:

{\sl (C1).} First select in the Mandelbrot set
a finite number
of truncated secondary limbs.
We require the hybrid classes of all quadratic-like renormalizations
 $R^m f$ to be picked from these limbs. 

{\sl (C2)}. On the other hand, we also require the 
{\sl combinatorial type} $\kappa(R^m f)$
to be sufficiently high on all levels (depending on the 
a priori choice of limbs).

\smallskip\noindent The second condition can be improved by specifying
other combinatorial factors producing a big space (see the subsection
with Lemma 11 and Lemma 16). 

\bigskip

\centerline{\bf \S 3. Geometric moduli.}\bigskip

Let us summarize the results of the section in two Theorems.
We say that a 
quadratic-like map $f: U'\rightarrow U$ has a definite modulus
if mod($U\backslash U')\geq\bar\mu>0$ (with an a priori selected quantifier 
$\bar\mu$).

\proclaim Theorem I. Let $f$ be a polynomial-like map with a definite modulus
whose internal class is selected from a given finite family of truncated
secondary limbs. Let $n(k)$ count the  levels  of non-central returns
(preceding the next quadratic-like level). Then the principle moduli 
$\mu_{n(k)+1}={\rm mod} (V^{n(k)}\backslash V^{n(k)+1})$
 grow with $k$ at uniformly linear rate.

\proclaim Theorem II. Let $f\in {\cal F}$.  Then all its quadratic-like
 renormalizations $R^nf$ have definite moduli. 

A compact set $K\subset {\bar C}$ is called {\sl removable} if given a 
neighborhood $U\supset K$, any
conformal embedding $\phi: U\backslash K\rightarrow \bar{\bf C}$ allows
the conformal continuation across $K$ (see [AB]). 
A simple condition  for
removability is the following. \smallskip

{\sl Assume that for any point $z\in K$
there is a nest of disjoint annuli $A_i\subset \bar{\bf C}\backslash K$
with definite moduli (mod($A_i)>\delta>0$) shrinking to $z$. Then
$K$ is removable.} \smallskip

Removable sets have zero Lebesgue measure. Now  Theorem II
immediately implies.

\proclaim Corollary IIa. Given an $f\in {\cal F}$,
its critical set $\omega(c)$ is a removable Cantor set.

By [HJ], [J] the a priori bounds also imply
the following (see the argument in \S 3). 

\proclaim Corollary IIb. The Julia set
   $J(f)$ of a map $f\in {\cal F}$ is locally connected.

According to McMullen [McM],  an infinitely q-renormalizable quadratic polynomial 
$f$ is called  {\sl robust} if for arbitrary high level $m$
there exists an annulus in ${\bf C}\backslash \omega(c)$
with definite modulus which is homotopic rel $\omega(c)$ to a Jordan curve
enclosing $J(R^m f)$
 but not enclosing any point of $\omega(c)\backslash J(R^mf)$.

\proclaim Corollary IIc. Any $f\in {\cal F}$ is
  robust. 

By [McM], robust quadratic polynomials have no 
invariant measurable line fields on the Julia set. 
Absence of invariant line fields  for a quadratic polynomial 
$f: z\mapsto z^2+c_0$ is equivalent to the property
that its topological class has empty interior  [MSS].
Theorem III below will  show that these topological classes are actually
single points for $f\in {\cal F}$.

\medskip\noindent{\bf Outline for Theorem II.} First we show that if a
quadratic-like map $f$ satisfying (C1) has a definite modulus then
the first  annulus of the principle nest also
has  a definite modulus. However the bound for this modulus is
certainly smaller than the a priori bound for $f$. 
To compensate this
loss, we go through the cascade of generalized renormalizations, and
observe (according to Theorem I)
 a linear growth of the principle moduli. So if we proceed for long
enough (assumption (C2)),
 we will arrive at the next quadratic-like level with a definite modulus
controlled by the same quantifier $\bar\mu$. Then we start over again. 

Most of this section is occupied with the proof of Theorem I.

\medskip\noindent{\bf Initial geometry.}

\proclaim Lemma 1. If the annulus $A$ has a definite modulus
 then the starting configuration
$(U,\gamma_i)$ of external rays has a bounded geometry.

\noindent{\bf Proof.} Indeed,  the map $f$ can be conjugate to a polynomial $g$
by a qc map with a bounded dilatation, where $g$ belongs to the finite set
of selected limbs. Let $g$ vary within one of these limbs. Then the 
finite intervals
of the external rays vary continuously with $g$. 

Since the truncated limbs don't touch the main cardioid, the absolute value of
the multiplier $\lambda$ of the $\alpha$-fixed point of $g$ is bounded
away from 1.
 Hence the fundamental annulus around this point
has a definite modulus. So the external rays landing at $\alpha$ will
meet this annulus on some definite distance from the Julia set.
Outside of the annulus they have a bounded geometry by
the previous argument. Near the fixed point the geometry is bounded
by a local consideration.         \QED 

Set $A^n=V^{n-1}\backslash V^n$.

\noindent{\bf Lemma 2.} The annulus $A^1=V^0\ssm V^1$
has a definite modulus (depending on the modulus  of $U\backslash U'$
only).

\noindent{\bf Proof.} Let us go back to the proof of Lemma 0.
Because of Assumption (C1), $V^0$ is well inside $\Omega^0$. As the puzzle-pieces
 $Y_j$
are obtained by pulling $V^0$ back by univalent iterates of $g$.
they are well inside $\Omega^0$ as
well. Finally, as $G^l$ in (0) is two-to-one branched covering,
$V^1$ is well inside of $Q$. 
   \QED 

\medskip\noindent{\bf A priori bounds.}

\proclaim Lemma 3. Let $i(1),...i(l)=0$ be the itinerary of a 
  puzzle-piece $V^{n+1}_j$ through the puzzle-pieces $V_i^n$ by 
   iterates of $g_n$
   until the first return back to $V^n$. Then
 $${\rm mod} (V^n\backslash V^{n+1}_j)\geq {1\over 2} 
     \sum_{k=1}^l{\rm mod} (V^{n-1}\backslash V^n_{i(k)}).$$

\noindent{\bf Proof.} The Gr\"{o}tcsz inequality.  \QED
 
 Let $D$ be a puzzle-piece which we call an "island" (compare below).
Let $W_i,\; i\in I$,  be a finite family of disjoint 
 puzzle-pieces  containing a critical puzzle-piece $W_0$.
We will freely identify the label set $I$ with the family itself.
For $W_i\subset D$ let 
$$R_i\equiv R_i(I, D)\subset D\backslash\bigcup _{j\in I}W_j$$ be
an annulus of maximal modulus enclosing $W_i$  but not enclosing other
puzzle-pieces of the family $I$. Such an annulus exists by the Montel
Theorem. We will briefly call it the {\sl maximal annulus} enclosing
$W_i$ in $D$ (rel the family $I$). 

Let us now define
{\sl the asymmetric modulus of the group $I$ in $D$ } as
$$\sigma(I|D)=\sum_{i\in I} \epsilon_i
  {\rm mod} (R_i), \eqno(1)$$
where the weight $\epsilon_i$ is equal to 1 for the critical puzzle-piece
and 1/2  for all others (if $D$ is a non-critical island then all weights are
actually 1/2).
This parameter  for a group of two puzzle-pieces was
suggested by Jeremy Kahn as a
 complex analogue of the asymmetric Poincar\'{e}
length [L2]. 

Let us now specify $D=V^{n-1}$, and $I$ to be a finite group of at least two 
puzzle-pieces $V^n_i$ of level $n$ containing the critical one.
Then set $\sigma_n (I)\equiv\sigma(I|V^{n-1})$ and
$$\sigma_n=\min_I\sigma_n(I), \eqno(2)$$
where $I$ runs over all groups of puzzle-pieces just specified.

Let us use a special notation for the {\sl principle moduli}
$$\mu_n={\rm mod} (V^{n-1}\backslash V^n). \eqno(3)$$
The $\mu_n$ and $\sigma_n$ are the principle geometric parameters 
of the renormalized maps $g_n$.

Our goal is to show that the asymmetric moduli  
 monotonically and linearly grow with $n$.
Let us fix a level $n\geq N$, denote $V^{n-1}=\Delta,\; V_i=V^n_i$,
$g=g_n$, and mark the
objects of the next level $n+1$ with prime. 

Let $I'$ be a finite family of puzzle-pieces $V'_i$.
Let us organize them in {\sl isles}
in the following way. Take two non-symmetric puzzle-pieces $V_i'$ and $V_j'$
and push them forward by iterates of $g$ through the puzzle-pieces $V_k$ of
the previous level. Find the first moment $t$ when they are separated by
those puzzle-pieces, that is, such that
$g^m V_i' $ and $g^m V'_j$ belong to the
same piece $V_{k(m)}$ for $m=0,...,t-1,$, while $g^t V_i' $ and $g^t V'_j$
land at different pieces. (In other words, the itineraries of  
$V_i'$ and $V_j'$ coincide until moment $t-1$). Then let us produce an
island $D$ by pulling $V_{k(t-1)}$ back by the corresponding inverse 
branch of $g^{t-1}$. Let $\phi_D = g^t : D\rightarrow\Delta$. This map
is either a double covering or a biholomorphic isomorphism
 depending on whether $D$ is critical or not. 

The family ${\cal D}={\cal D}(I')$ of isles form a lattice with respect 
to inclusion. Let 
{\sl depth} $: {\cal D}\rightarrow {\bf N}$ be the minimal strictly
 monotone function on this lattice, assigning to the biggest
island  $V_0\equiv \Delta'$ depth 0.

Let us now consider the asymmetric moduli $\sigma(I|D)$ as a function on
the family ${\cal D}$ of isles. This function is clearly monotone:
$$\sigma(I|D)\geq\sigma(I|D_1)\quad {\rm if} \quad D\supset D_1, \eqno (4)$$
 and superadditive:
$$\sigma(I|D)\geq\sigma(I|D_1)+\sigma(I|D_2),   \eqno  (5)$$
provided $D_i$ are disjoint subisles in $D$.

We call a  puzzle-piece $V_j'\subset D$ {\sl pre-critical rel} $D$ 
if $\phi_D(V_j')=V_0$. If $D=\Delta'$ is the trivial island, we skip "rel".
There are at most two pre-critical pieces in any $D$.
If there are actually two of them, then they are non-critical and
symmetric with respect to $c$. 

Let $D$ be a deepest island of family ${\cal D(I')}$, and
 $V_j',\; j\in J$, be the
group of puzzle-pieces contained in $D$, that is $J=I'|D$.
Let $i(j)$ is defined for $j\in J$ by the property
 $\phi_D(V_j')\subset V_{i(j)}$, and $I=\{i(j) : j\in J\}$.

\proclaim Lemma 4.
  Under the circumstances just described the following estimate holds:
$$ \sigma(J|D)\geq 
{1\over 2}\left( (|J|-s)\mu +  s\; {\rm mod}(R_0)+
 \sum_{i\in I, i\not=0}{\rm mod}(R_i)\right), 
\eqno (6)$$ 
where $s=\#\{j: i(j)=0\}$ is the number of pre-critical pieces rel $D$.

\noindent{\bf Proof.}  As $D$ is the deepest
island, each puzzle-piece $V_i,\; i\in I,$ contains a single puzzle-piece
$\phi_D V_j$ (though there might be two symmetric puzzle-pieces 
in $J$ with $\phi_D V_j=\phi_D V_k$). Let $R_i\subset \Delta$
 denote an annulus of maximal modulus enclosing $V_i$ rel $I$, and
 let $T_j\subset D$ be an annulus of maximal modulus
enclosing $V'_j$ rel family $J$.
Let $\delta_{st}$ denote the Kronecker symbol.
Fix a $j\in J$ and let $i=i(j)$.
Let us consider now two cases:

\smallskip
(i) Let $V_j'$ be non-critical.  Then
    $${\rm mod }(T'_j)\geq {\rm mod}(R_i)+\delta_{0i}\;\mu.  \eqno(7)$$
    To see that, observe that mod$(V_i\backslash \phi_D V'_j)$ is at least
    $\mu$, provided $i\not=0$. 
    Observe also that the pull-back of the topological disc
    $Q_i=R_i\cup V_i$ to $D$ is
    univalent. Indeed, if $\phi_D$ were a double covering then the island $D$
    would be critical, and hence  would contain the critical
    puzzle-piece $V_0'$. It follows that
    $Q_i$ does not contain the critical value of $\phi_D$.

\smallskip
(ii) Let $V'_j=V'_0$ is critical. Then
   $${\rm mod} (T_0')\geq {1\over 2}({\rm mod}(R_i)
    +\delta_{0i}\;\mu). \eqno (8)$$

Summing up the estimates (7) and (8) with the weights 1/2 and 1
correspondingly over the family $J$, we obtain the desired estimate.  \QED

\proclaim Corollary 5. For any island $D$ of the family ${\cal D}(I')$
 the following estimates hold:
$$\sigma(I'|D)\geq {1\over 2}\mu\quad {\rm and}
\quad \sigma(I'|D)\geq\sigma.$$

\noindent{\bf Proof.} By  monotonicity (4), it is enough to check the case
of a deepest island $D$. Let us use the notations of the previous lemma.
Observe first that the family $I=\{i(j): j\in J\}$ contains at least two
puzzle pieces. Indeed, the only case when $|I|<|J|$ can happen is when
$\phi_D$ is a double covering, and there are two symmetric puzzle-pieces in 
the family $J$. But then this family must also contain the critical piece
$V'_0$, and hence $|I|>2$.

As $\mu>{\rm mod}(R_0)$, $|J|\geq 2$ and $|I|\geq 2$,
 the right-hand side in (6) is bounded from
 below by
$$ {1\over 2}\left( |J|\; {\rm mod}(R_0) + \sum_{i\in I, i\not=0}
 {\rm mod}(R_i)\right)\geq \sigma(I)>\sigma,  \eqno (9)$$ 
  \QED      

Let us decompose $g_n: V^n\rightarrow V^{n-1}$ 
as $h_n\circ\Phi$ where $\Phi$ is purely
quadratic, while $h_n$ is univalent. The {\sl non-linearity} or
{\sl distortion} of $h_n$ is defined as
$$\max_{z,\zeta\in V^n} \log\left|{Dh_n(z)\over Dh_n(\zeta)}\right|,$$
and measures how far $g_n$ is from being purely quadratic.

\proclaim Corollary 6 (a priori bounds). 
 The asymmetric moduli $\sigma_n$ grow monotonically
  and hence stay away
  from 0 on all levels (until the next quadratic-like level).
   The basic moduli $\mu_n$ stay away
   from 0 everywhere except for  tails of long cascades
   of central returns. Moreover, the non-critical puzzle-pieces $V^n_i$
   are also well inside $V^{n-1}$ except for pre-critical pieces on the
    levels which immediately follow the long cascades of central returns.
 The distortion of $h_n$ is uniformly bounded on all levels.

\noindent{\bf Proof.} On the first non-degenerate level $N+1$ we have a
definite principle modulus by Lemma 2. Hence by the previous Corollary we have 
a definite value of $\sigma$ on the next level which then begins to grow
monotonically. So, it stays definite on all levels until the next 
quadratic-like one. By Lemma 3, the basic moduli stay definite
as well, except for tails of  long cascades
   of central returns.  The next statement also follows from Lemma 3. 

To check the last statement, it is enough to observe that $h_n$ has a
Koebe space spread over $V^{n-2}$. Hence its distortion is controlled
by the principle scaling factor $\mu_{n-1}$. So we are OK outside the tails
of central cascades. But observe also that within the central cascade
we keep the same return map, just shrinking its domain.    \QED

\medskip\noindent{\bf Linear growth.}
Our goal is to prove that $\sigma'\geq\sigma+a$ with
a definite $a>0$ at least on every other level except for the tails
of central cascades.
Corollary 6 shows the reason why these tails play a special
role. The  growth rate of $\sigma$ definitely slows down in the
tails. So let us assume that the level
 $n-1$ is not there, so that the principle modulus $\mu$ is definitely
positive.

\noindent{\bf Corollary 7.} If a deepest island $D$ contains
  at least three puzzle-pieces $V'_j,\; j\in J$,  then
  $$\sigma(J|D)\geq\sigma(I)+{1\over 2}\mu.$$

\noindent{\bf Proof.}  Let us in (6) split off
$(1/2)\mu$ and estimate all other $\mu$'s by mod$(R_0)$. 
This estimates the right-hand side by
$${1\over 2}\mu+{|J|-1\over 2}{\rm mod}(R_0)+
{1\over 2}\sum_{i\in I, i\not=0}{\rm mod}(R_i), $$
which immediately yields what is claimed.     \QED

Let us now consider the case when the island $D$ contains
only two puzzle puzzle-pieces. In order to treat it, we need some preparation
in geometric function theory.

\medskip\noindent{\bf Moduli defect, capacity and eccentricity.} 
Let $D$ be a topological disk, $\Gamma=\partial D$,
 $a\in D$, and $\psi: (D,a)\rightarrow
({\bf D}_r,0)$ be the Riemann map onto a round disk of radius $r$
with $\psi'(a)=1$.
Then $r\equiv r_a(\Gamma)$ is called the {\sl conformal radius} of $\Gamma$ 
about $a$. The {\sl capacity} of $\Gamma$ rel $a$ is defined as
$${\rm cap}_a(\Gamma)=\log r_a(\Gamma).$$

\proclaim Lemma 8. Let $D_0\supset D_1\supset K$, where $D_i$ are topological
disks and $K$ is a connected compact. Assume that 
the hyperbolic diameter of $K$ in $D_0$ and the hyperbolic 
dist$(K,\partial D_1)$ are both bounded by a $Q$.
  Then there is an $\alpha(Q)>0$
such that
$${\rm mod }(D_1\backslash K)\leq {\rm mod} (D_0\backslash K)-\alpha(Q).  $$

\noindent{\bf Proof.}  Let us take a point $z\in \partial D_1$ whose hyperbolic
distance to $K$ is at most $Q$. Then there is an annulus  of a definite
modulus contained in $D_0$ and enclosing both $K$ and $z$. 

Let us uniformize $D_0\backslash K$ by a round annulus
 $A_r=\{\zeta : r<|\zeta|<1\},$
and let $\tilde z$ correspond to $z$ under this uniformization. Then 
$\tilde z$ stays a definite Euclidian distance $d$ from the unit circle.

If $R\subset A_r$  is any annulus enclosing the inner boundary of $A_r$
but not enclosing $\tilde z$ then by the normality argument
mod$(R)<{\rm mod}(A_r)-	\alpha_r(d)$ with an $\alpha_r(d)>0$.
(Actually, the extremal annulus is just $A_r$ slit along the radius from
$\tilde z$ to the unit circle).

We have to check that $\alpha_r(d)$ is not vanishing as $r\to 0$.
Let us fix an outer boundary $\Gamma$ of $B$ (the unit circle + the slit 
in the extremal case). We may certainly
assume that the inner boundary coincides with the $r$-circle.
Then the defect mod$(R)-\log(1/r)$ monotonically increases to the
 cap$_0(\Gamma)$. By normality this capacity is   bounded above
 by an $-\alpha(d)<0$, and we are done.      \QED

Let $A$ be a standard cylinder of finite modulus, $K\subset A$. Let us define
the  mod($K$) as the modulus of the smallest concentric sub-cylinder 
$A'\subset A$ containing $K$ (see Figure 3).

\proclaim Lemma 9 (Definite Gr\"{o}tcsz Inequality). Let $A_1$ and $A_2$
be homotopically non-trivial disjoint topological annuli in $A$.
Let $K$ be the set of points in their complement which are separated 
by $A_1\cup A_2$ from
the boundary of $A$. Then there is a function $\beta(x)>0$
($x>0$) such that
$${\rm mod }(A)\geq {\rm mod}(A_1)+{\rm mod}(A_2)+\beta({\rm mod}(K)).$$

\noindent {\bf Proof.} For a given cylinder this follows from the
usual Gr\"{o}tcsz Inequality and the normality argument. 
Let us fix a $K$, and let  mod$(A)\to\infty$. We can assume that $A_i$
are lower and upper components of $A\backslash K$ correspondingly.
Then the modulus defect
$${\rm mod }(A)-{\rm mod}(A_1)-{\rm mod}(A_2)$$
 decreases by the usual Gr\"{o}tcsz inequality. At the limit the cylinder
becomes the punctured plane, and the modulus defect converges to
-(cap$_0(K)+$cap$_{\infty}(K)$). It follows from the area estimates that
this sum of capacities is negative, unless $K$ is a circle centered at the
origin. This estimates depends only on mod$(K)$ by normality.  \QED

Let $d_a(\Gamma)$ and $\rho_a(\Gamma)$ be the Euclidian radii of the 
inscribed and circumscribed circles about  $\Gamma$ 
centered at $a$. Then let us define the eccentricity of $\Gamma$ about $a$ as
$$e_a(\Gamma)=\log{\rho_a(\Gamma)\over d_a(\Gamma)}.$$ By Koebe and
Schwarz,
$$e_a(\Gamma)=-({\rm cap}_a(\Gamma)+{\rm cap}_{\infty})(\Gamma)+O(1),$$ 
with $O(1)\leq 2\log4$.

\proclaim Lemma 10. Under the circumstance of Lemma 9 assume also that 
the annulus
$A\subset {\bf C}\backslash \{a\}$ is homotopically non-trivially embedded
in the punctured plane, and mod$(A_i)\geq\alpha>0$. If $e_a(K)$ is big
then mod$(K|A)$ is big as well.

\noindent{\bf Proof.} Let us  consider the uniformization $\phi:\bar A
\rightarrow A$ of $A$ by a round annulus. If mod$(K|A)$ is bounded then
$\bar K$ is well inside of $\bar A$. Then by the normality argument
$K$ must have a bounded eccentricity about $a$.   \QED   

\medskip
\noindent{\bf The case of two puzzle-pieces.} Let us now go back to the
estimates of asymmetric moduli.  Suppose we have a deepest island $D$ containing
two puzzle-pieces $V_j^{n+1},\; j\in J$. Let $\phi\equiv \phi_D$ and let
$\phi V^{n+1}_j\subset V_i^n$ with $i=i(j)$. Let us split the
argument into several cases.

\smallskip {\sl Case (i)}.
   There is a non-critical puzzle-piece $V^n_i,\; i\in I$, 
   which stays on a bounded 
  Poincar\'{e} distance in $V^{n-1}$ (controlled by a given big quantifier $Q$)
 from the critical point. Then by Lemma 8 
$$\mu_n\geq {\rm mod}(R_0) + \alpha  \eqno (10)$$ with a definite
$\alpha=\alpha(Q)>0$. But observe that when we passed from Lemma 4
to Corollary 5 we estimated $\mu$ by mod$(R_0)$. Using the better
estimate (10), we obtain a definite increase of $\sigma$. 

\smallskip {\sl Case (ii).}
   Let each non-critical puzzle-piece $V^n_i,\; i\in I$,
  stay hyperbolically far away from the critical point.
Then $V^n_0$ may not  belong
to any non-trivial island together with some non-critical piece
  $V^n_i,\; i\in I$. Indeed, it follows from Corollary 6 that
any non-trivial island is well inside of $V^{n-1}$.

 {\sl Assume first that both $V^n_i$ are non-critical}. Then $\sigma(J|D)$
is estimated by $\sigma_n(\tilde I)$ where the family $\tilde I$  consists
of $V^n_i$ and the central puzzle-piece $V^n_0$. If no two of these
puzzle-pieces belong to the same non-trivial island, then by Corollary 7
$\sigma(\tilde I)\geq \sigma_{n-1}+a$ with a definite $a>0$.

 Otherwise the puzzle-pieces $V^n_i, \; i\in I,$ belong to an island $W$. 
Since $W$ is well inside of $V^{n-1}$, it stays on the big Poincar\'{e}
distance from the critical point. Hence mod$(R_0)\approx\mu$
(this sign means the equality up to a small constant controlled by the
quantifier $Q$, while the sign $\succ$ below means the inequality up to
a small error), and
$$\sigma(\tilde I)\geq \sigma(I|Q)+{\rm mod}(R_0)\succ
   \sigma^{n-1}+\mu.$$
So we have  gained some extra growth, and can pass to the next case.

\smallskip\noindent
 {\bf Fibonacci returns.}
  {\sl Let one of the puzzle-pieces $V_i^n$ be critical}.
  So we have the family $I^n$ of two puzzle-pieces $V_0^n$ and $V_1^n$.
  Remember that we also assume that {\sl the hyperbolic distance
  between these pieces is big}. Hence, $V^{n-1}$ is the
  only island containing both of them, so that $g_{n-1} V_0^n$ and 
 $g_{n-1} V_1^n$
  belong to different puzzle-pieces of level $n-1$.  For the same reason 
   we can assume that one of these puzzle-pieces is critical. Denote
  them by $V^{n-1}_0$ and $V^{n-1}_1$. 
  Then one of the following  two possibilities on level $n-2$ can occur:

\smallskip\noindent
1) {\sl Fibinacci return} 
   when $g_{n-1} V^n_0\subset V^{n-1}_1$ and
   $g_{n-1}V^n_1=V^{n-1}_0$ (see Figure 4);

\smallskip\noindent 2) {\sl Central return}  when
    $g_{n-1} V^n_0= V^{n-1}_0$ and
   $g_{n-1}V^n_1\subset V^{n-1}_1$. 

  We can assume that one of these schemes occur on several previous levels
    $n-3, n-4,...$ as well (otherwise we gain an extra growth by the previous 
   considerations). To fix the idea, let us first consider
  the following particular case

\smallskip\noindent{\sl Fibonacci cascade.} Assume that  on both
   levels $n-1$ and $n-2$ the Fibonacci returns  occur. 
  Let us look more carefully at the estimates of Lemma 4. 
 In the Fibonacci case we just have:

$${\rm mod}(R^n_1)\geq {\rm mod} (R^{n-1}_0),  \eqno (11)$$
$${\rm mod} (R^n_0)\geq {1\over 2} {\rm mod} 
   (g_{n-1} V^n_0 | Q^{n-1}_1), \eqno (12)$$
where $Q^n_i=V^n_i\cup R^n_i$. Applying $g_{n-2}$ we see that
$${\rm mod} ( Q^{n-1}_1\backslash g_{n-1} V^n_0 )\geq
 {\rm mod} (Q^{n-3}_0\backslash V_0^{n-1}). \eqno (13)$$
But since $V^{n-2}_1$ is hyperbolically far away from the critical point,
$${\rm mod} (Q^{n-3}_0\backslash V_0^{n-1})\approx
 {\rm mod}(V^{n-3}_0\backslash V_0^{n-1}). \eqno (14) $$
By the Gr\"{o}tcsz Inequality there is an $a\geq 0$ such that
$${\rm mod}(V^{n-3}_0\backslash V_0^{n-1})=\mu_{n-1}+\mu_{n-2}+a.  \eqno (15)$$
Clearly
$$\mu_{n-1}\geq {\rm mod} (R_0^{n-1})
. \eqno (16) $$
Furthermore, let $P_1^{n-1}$ be the pull-back of $Q_0^{n-2}$ by $g_{n-2}$.
Since $\partial(P_1^{n-1})$ is hyperbolically far away from $V_1^{n-1}$,
we have:
$$\mu_{n-2}\geq {\rm mod}(R_0^{n-2})={\rm mod} (V_1^{n-1}|P_1^{n-1})
      \approx {\rm mod} ( R_1^{n-1}). \eqno(17) $$
Combining estimates (12)-(17) we get
$${\rm mod}(R_0^n)\succ {1\over 2} ({\rm mod}(R_0^{n-1})+{\rm mod}(R_1^{n-1})+a)     .  \eqno (18)$$

We see from (11) and (18) that the only thing to check  that
the constant  $a$ in (15) is definitely positive.
Assume that this is not the case. Set $\Gamma_n=\partial V^n$.
Then by the Definite Gr\"{o}tscz Inequality, the mod$(\Gamma_{n-2})$
in the annulus $A=V^{n-3}\backslash V^{n-1}$ is very small. 
Since $\Gamma_{n-2}$ is well inside of $A$, we conclude by the Koebe
Distortion Theorem that $\Gamma_{n-2}$ is contained in a narrow
neighbourhood of a curve $\gamma$ with a bounded geometry. Moreover, this
curve has a definite eccentricity around the critical point.

On the other hand, the puzzle-piece $V^{n-1}_1$ is hyperbolically far away
from the critical point. Hence it must be located  Euclidianly  very close 
 to $\Gamma_{n-2}$ (relatively the Euclidian 
distance to the critical point). Hence the critical value $g_{n-1} c$
is also extremely close to $\Gamma_{n-2}$.

By Corollary 6,
$g_{n-1}$ is  a quadratic map up to a bounded distortion. Hence the
curve $\Gamma_{n-1}$ which is the pull-back of $\Gamma_{n-2}$ by $g_{n-1}$
must have a huge eccentricity around the critical point. By Lemmas 10 and 9
it will contribute towards the definite extra constant on the $(n+1)$st level.

\smallskip\noindent {\bf Remark.} The actual shape of a deep level puzzle-piece
  for the Fibonacci cascade is shown on Figure 5. There is a good reason 
 why it resembles the filled-in Julia set for $z\mapsto z^2-1$
  (see [L3]). As the geodesic in $V^{n-1}_0$ joining 
   the puzzle-pieces $V^n_0$ and $V^n_1$
   goes  through the  pinched region,
 the Poincar\'{e} distance
  between  these puzzle-pieces is, in fact, big.
                                                 
\smallskip\noindent{\sl General case.}
Let us now allow the central returns
  along with the Fibonacci ones. Suppose we have
a cascade of central returns on $N-1$ subsequent levels
$ V^{m}\supset...\supset V^{m+N-2}\equiv V^{n-2}$, preceded by
 the Fibonacci return on level $m-1$. So
$g_{m+1} c\in V^{m+N-1}_0$, while $g_m c\in V^m_1.$ By our convention,
   this cascade is not too  long, so that we have a definite
  space in between any two levels. 

Let us now pass from the island $D\subset V^n\equiv V^{m+N}$ all way up
the cascade to the level $m-1$, that is, consider the map
 $$G=g_m\circ g_{m+1}^{N-1}\circ \phi_D : D\rightarrow V^{m-1}. \eqno (19)$$
Then $S\equiv G\; V^{m+N+1}\subset V^m$. Now we again should split the analysis
depending on where the puzzle-piece $ V^{m+N+1}$ lands. 
Let us start with the most
interesting case when  it lands at the deepest possible level.

\smallskip {\sl Subcase (a).} Let  $S=V^{m+N}$. Pulling the annuli 
$R^m_0$ and $R^m_1$  back by $G$ to $D$, we get the following estimates:

$$ {\rm mod} (R_0^{m+N+1})\succ {1\over 2} {\rm mod}(V^{m-1}\backslash S)
   ={1\over 2} {\rm mod}(V^{m-1}\backslash V^{m+N}), \eqno (20)$$

$${\rm mod} (R^{m+N+1}_1)\geq {1\over 2^N} 
  ({\rm mod}(R_1^m)+ {\rm mod}(R_0^m)).
   \eqno (21)$$ 
By the Gr\"{o}tcsz inequality there is an $a\geq 0$ such that
$${\rm mod}(V^{m-1}\backslash V^{m+N})\geq  {\rm mod}(A^m) +
    {\rm mod}(V^m\backslash V^{m+N})+a\geq$$
    $${\rm mod} (R_0^m)+\sum_{k=m+N}^{m+1}{\rm mod}(A^k)+a\geq
    {\rm mod} (R_0^m)+\sum_{k=0}^{N-1}{1\over 2^k}{\rm mod}(R_0^{m+1})+a. 
   \eqno (22)$$
Since ${\rm mod}(R_0^{m})\approx {\rm mod}(R_1^{m+1})$,
the above  estimates imply

$$2\,\sigma(I^{m+N+1}|D)=2\, R_0^{m+N+1}+ R^{m+N+1}_1\succ$$
$$  {1\over 2^N} {\rm mod}(R_1^m)+
   {1\over 2^{N-1}}{\rm mod} (R_0^m)+
   (1-{1\over 2^N})  {\rm mod}(R_1^{m+1})
 +(2-{1\over 2^{N-1}})   {\rm mod}(R_0^{m+1})+a\approx$$
  $$\approx {1\over 2^{N-1}} \sigma(I^m)+
  (2-{1\over 2^{N-1}})\sigma(I^{m+1})
    +a .$$

We see that if the curve $\Gamma^m$ has a definite modulus in the annulus
$V^{m-1}\backslash V^{m+N}$ then we have a definite growth of $\sigma$.
Otherwise arguing as in the case of the Fibonacci cascade we 
conclude that  the curve $\Gamma^{k}$ has a big eccentricity around
 the  puzzle-piece $V^{k+1}_1,\; k=m,...,m+N-1$.

Let us now go one central cascade up  to the level $V^{m-1}$ (until the 
Fibonacci level). If this cascade is not too long,  
then by the above considerations 
we  either  have a definite growth of $\sigma$ 
within this cascade, or  $\Gamma^{m-1}$ has a big eccentricity about
$V^m_1$. But then $\Gamma^m$ has a big modulus in $V^{m-1}\backslash V^{m+N}$,
and we are done.

Finally, if $V^{m-1}$ is in the tail of a long central cascades 
then $\Gamma^m$ has {\sl always}  a big eccentricity about the critical
point (see the next subsection). If we actually have a central return on
level $m$ (so that $N\geq 2$), then $\Gamma^{m+1}$ has a big eccentricity
around $c$ as well. But this curve is for sure well inside 
$V^m\backslash V^{m+N}$. So we can use it instead of $\Gamma^m$
to contribute to the definite $a$ in estimate (22).

If a non-central return on level $m$ occurs (that is, $N=1$),
 then we don't see a definite
growth for $\sigma_{m+N+1}$ but we gain it one level down.

\smallskip{\sl Subcase (b).}  Assume now that
   $S\subset V^m\backslash V^{m+N}.$  Let us consider the
  Markov family of puzzle-pieces $W^k_i,\; k=m+1,...,m+N,$ the pull-backs
  of pieces $V^{m+1}_i\equiv W^{m+1}_i$ to the annuli $A^k$. 
  Let $S\subset W\equiv W^k_i$. Then
$${\rm mod }(W\backslash S)\geq {\rm mod}(V^m\backslash V^{m+N}),$$
  and  we have
  $${\rm mod} (V^{m-1}\backslash S)\geq
 {1\over 2} ({\rm mod}(A^m)+{\rm mod}(W\backslash S)+
  {\rm mod}(V^m\backslash W))\succ$$
$${1\over 2} ({\rm mod}(R^m_0)+{\rm mod}(V^m\backslash V^{m+N})+a$$
where $a>0$ is definite, unless $V^{m-1}$ is in the tail of
 a long central cascade. But then argue as in Subcase (a).
Theorem A is proved.

\medskip\noindent {\bf Other factors yielding big space.}
Theorem A ensures that mod$(Rf)$ is sufficiently high if the type 
is sufficiently high, that is, there are 
sufficiently many non-central levels. However, there are other
combinatorial factors which imply big mod$(Rf)$ as well.
For example, if the return time of some $V^{n+1}_j$ back to
$V^n$  under iterates of $g_n$ is high, then Lemma 3 implies big space.

Sometimes long central cascades imply big space as well. 
Let us  consider such a cascade 
$$ V^{m}\supset...\supset V^{m+N-1},$$ 
where $gc\equiv g_{m+1} c\in V^{m+N-1}\backslash V^{m+N}$. 
The quadratic-like map $ g: V^{m+1}\rightarrow V^m$ can be viewed as
a small perturbation of a quadratic-like map $G$ with a definite modulus
and with non-escaping
critical point. Let $c\in \partial M$ be the internal class of $G$.

\proclaim Lemma 11. Under the above circumstances let us assume that
$z\mapsto z^2+c$ does not have neither parabolic points nor Siegel disks.
If $g$ is sufficiently close to $G$ (depending on $c$ and a priori bounds)
then mod$(A^{m+N+3})$ is big.

\noindent{\bf Proof.} The above assumptions mean that the Julia set
$J(G)$ has empty
interior.  If $g$ is sufficiently close to $G$ 
then  $\Gamma^{m+N-1}=\partial V^{m+N-1}$
is close in the Hausdorff metric to $J(G)$.
Hence  $\Gamma^{m+N-1}$ has a big eccentricity with respect to any point
$z\in V^{m+N-1}$.

As $g_m$ are quadratic maps up to bounded distortion, the curves
$\Gamma_{m+N}$, $\Gamma_{m+N+1}$ and $\Gamma_{m+N+2}$ also have big eccentricity
with respect to any enclosed point. Moreover, there is a definite
space in between these two curves. Hence mod$(V^{m+N}\backslash V^{m+N+2})$
is big. This implies that mod$(A^{m+N+3})$ is big as well.

Indeed, if central return on level $m+N$ occurs then 
straightening the quadratic-like map 
$g_{m+N+1}: V^{m+N+1}\rightarrow V^{m+N}$ by a qc map we conclude that
$${\rm mod}( V^{m+N}\backslash  V^{m+N+2})\asymp {\rm mod}(A^{m+N+1}).$$
Hence $A^{m+N+1}$ has a big modulus.


So we can assume that non-central returns occur on levels $m+N$  and
$m+N+1$. Let us show that then mod$(A^{m+N+3})$ is big.
Let $\psi^{\circ s} c\in V_j^{m+N+2}$. Then it is easy to see that
$${\rm  mod}(A^{m+N+3})\geq 
 {1\over 2}{\rm mod} (V^{m+N+1}\backslash \psi^{\circ s}( V^{m+N+3})).$$
Let now $t$ be the return time of $V^{m+N+2}_j$ back to $V^{m+N+1}$ under
iterates of $g_{m+N+1}$. Under this iterate  $\psi^{\circ s} (V^{m+N+3})$
 is mapped
onto $V^{m+N+2}$, and we conclude that
$${\rm  mod}(A^{m+N+3})\geq 
 {1\over 4}{\rm mod} (V^{m+N}\backslash V^{m+N+2}),$$
which is big.  \QED

\smallskip\noindent{\bf Remark.}
In the real case we will give a complete description of the combinatorial
factors producing big space (see Lemma 16).


\medskip\noindent{\bf Proof of Corollary IIb. (local connectivity
of the Julia sets)}.  I learned the
following argument from J.~Kahn and C. McMullen.
 It follows from Theorem II  that the
renormalized Julia sets $J(R^m f)$ shrink down to the critical point.
Let us take an $\epsilon>0$, and find an $m$ such that $J(R^m f)$ is
contained in the $\epsilon$-neighborhood of the critical point.

Let $\alpha_m$ denote
the dividing fixed point of the Julia set $J(R^m f)$, and $\alpha_m'$
denote the symmetric point. Let us consider a topological disk
bounded by an equipotential level, and cut it by the
external rays landing at $\alpha_0,...,\alpha_{m-1}$ into the puzzle-pieces
$P^{m,1}_j$ (as Yoccoz did in the finitely q-renormalizable case).
 Let us then pull these puzzle-pieces back in the usual way,
and use the notation $P^{m,l}(a)$ for the puzzle-piece of level $l$
containing a point $a$.

Consider the nest $P^{m,1}(c)\supset P^{m,1}(c)\supset...$ of the critical
puzzle-pieces. This nest shrinks down  to the Julia set $J(R^mf)$. 
Hence there is a puzzle-piece $P^{m,l}(c)$
 contained in the $\epsilon$-neighborhood
of the critical point. As $J(f)\cap P^{m,l}(c)$ is clearly connected,
the Julia set $J(f)$ is locally connected at the critical point.

Let us  now prove local connectivity at any other point $z\in J(f)$.
Consider two cases.\smallskip

{\sl Case (i)}. Let the orbit of $z$ eventually land at all Julia sets
  $J(R^mf)$.  Take the first moment $k=k(m)$ such that $f^kz\in J(R^mf)$.
Let us show that the domain $U$
can be univalently pulled back along the orbit $z,...,f^k z$.
Let $U_m'\equiv V^{m,t(m)}$, $U_m\equiv V^{m,t(m)-1}$, 
$p$ be the return time of $c$ back to $U_m$, and
$${\bf Q}_m\equiv \bigcup_{t=1}^p f^t U_m' \eqno (23).$$
Let us find the smallest
natural number $l$  such that 
$f^l z\in {\bf Q}_m,$ and  moreover let
 $f^l z\in f^s U_m, \; 1\leq s\leq p$. Then $f^{l-1} z$ belongs to the
domain $\Omega$ which is $c$-symmetric to $f^{s-1} U_m$.
 As $\Omega$ is disjoint
from  ${\bf Q}_m\supset \omega(c)$, there is a single-valued branch
$f^{-l}: \Omega\rightarrow Z\ni z$. On the other hand, clearly there is a 
single-valued branch $f^{-(s+1)} : U_m\rightarrow\Omega$. Hence there is
a single-valued branch $f^{-k} : U_m\rightarrow Z$ as it was claimed.

Because of the a priori bounds, the Julia set $J(R^m f)$ is well inside
of $U_m$. Hence there is a puzzle piece $P^{m,l}(c)\supset J(R^mf)$ which is
well inside of $U_m$ as well. It follows from the Koebe Theorem that
its pull-back $Y\ni z$ has a bounded shape and hence a small diameter
(for sufficiently big $m$). As $Y\cap J(f)$ is connected, we are done.

\smallskip{\sl Case (ii)}. Assume that the orbit of $z$ never lands at
  $J(R^mf)$. Then it never lands at the forward orbit 
$J_m$ of $J(R^m f)$. Hence it accumulates on some point $a\not\in J_m$.
  But the puzzle-pieces
  $P^{m,l}(a)$ are disjoint from the critical set for sufficiently big $l$. 
Pulling them back
to $z$, we again obtain small pieces $Y\ni z$ containing a
 connected  part of the Julia set.           \QED

\bigskip\centerline{\bf \S 4. Pull-Back Argument.} \medskip

Any quadratic polynomial induces an equivalence relation on the
rational points of the circle ${\bf T}$ by identifying the external
arguments whose external rays land at the same point of the Julia set
(see Douady and Hubbard [DH1], [D2] and [H]).
Two quadratic-like maps are called {\sl combinatorially equivalent} if
they induce the same equivalence relation (the combinatorial classes are
clearly bigger than the topological ones). The combinatorial class of
quadratic polynomials is obtained by intersecting a
nest of parameter puzzle-pieces bounded by appropriate external rays
and equipotentials. The definition of combinatorially equivalent 
quadratic-like maps is straightforward.
 
Our goal is to prove the following result.

\proclaim Theorem III. Let  $f$ and $\tilde f$ be two quadratic-like maps
 of class ${\cal F}$. If these maps are combinatorially
  equivalent then they are quasi-conformally conjugate.

\proclaim Corollary. Any  quadratic polynomial $f: z\mapsto z^2+c$ of 
  class ${\cal F}$ is combinatorially rigid, so that MLC holds at $c$.

\noindent {\bf Proof of Corollary}. 
The well-known argument: combinatorial classes of quadratic polynomials 
are closed, while qc classes are either open or single points. So if
a combinatorial class coincides with a qc class, both must be single points.
MLC follows since the intersection of the Mandelbrot set with puzzle-pieces
is connected.  \QED

\medskip\noindent {\bf Strategy.} The method we use for proof of 
Theorem III is called ``the pull-back argument". The idea is to start
with a qc map respecting some dynamical data, and then pull it back 
so that it will respect some new data on each step. In the end 
it becomes (with some luck)  a qc conjugacy. This method originated in
Sullivan's work, and then was developed 
in several other works (see [K] and [Sw]). Our way is to pull back
through the cascade of generalized renormalizations. The linear
growth of moduli gives us enough dilatation control until the next
quadratic-like level, while complex a priori bounds
allow us to interpolate  and pass to the next level.

We will use tilde 
for marking the corresponding objects. Referring to a qc-map, we always mean
that it has a definite dilatation. All puzzle-pieces have a natural
boundary marking coming, e.g, from the uniformization of the basin at $\infty$
(we can always assume that we have started with a polynomial map). 
Let us call two configurations of puzzle-pieces
$W_i$ and $\tilde W_i$ qc pseudo-conjugate
if  there is a qc map between them respecting the boundary marking.  

 Let $\{V^{m,n}\}$ be the principle nest of
critical puzzle-pieces (see \S 1). We switch from the
 nest $V^{m,n},\; n=0,1,...,t(m),$ to the next nest
$V^{m+1,n},\; n=0,1,...,$
when the modulus $A^{m,t(m)+1}$ is bounded from the both sides (not only
from below). 

By Lemma 1  the starting
configurations $\{V^{m,0}_i\}$  and $\{\tilde V^{m,0}_i\}$ have 
bounded geometry, so there is a qc pseudo-conjugacy $h_m$ between them. 
It is possible to pull it back to the first non-degenerate level,
no matter how deep it is (The Initial Construction below).
Let us then
pull it back through the cascade of generalized renormalizations
(the Main Step below). Since the geometric moduli of these
maps linearly increase, the positions of their critical values 
are  localized with an exponentially high precision. 
It follows that the qc dilatation of the pseudo-conjugacy
on the next level can jump only by an exponentially
small amount. Hence we will arrive at the next quadratic-like level $m+1$
with a qc map $H_m$ with bounded dilatation.  

Finally, since the annuli  $A^{m,t(m)+1}$ and  $\tilde A^{m,t(m)+1}$
  have a definite moduli,
we can qc interpolate in between $H_m$ on their outer boundaries and
$h_{m+1}$ on the inner ones (keeping the map in the right homotopy class
mod the critical set). This gives us a qc pseudo-conjugacy between
the  nests of critical puzzle-pieces. Then it can be easily spread around 
to the whole critical set. Sullivan's pull-back argument completes the 
construction.

\medskip\noindent{\bf Main Step.} 
Let $g: \cup V_i\rightarrow\Delta$ and 
$\tilde g: \cup\tilde V_i\rightarrow \tilde \Delta$
 be two generalized polynomial-like maps. The objects on the next
 renormalization level will be marked with prime.
So $g': \cup V_j'\rightarrow \Delta'$ is the generalized renormalization
of $g$, $\Delta'\equiv V_0$. 
 Let $\mu$ the principle modulus of $g$.

\smallskip\noindent {\bf Remark.} We don't assume that the non-critical 
puzzle-pieces $V_i^n,\; i\not=0,$ are well inside $\Delta$, since this is
not the case on the levels which immediately follow long cascades of
central returns. We even allow the annuli $\Delta\backslash V^n_i,\; i\not=0$
to be degenerate which actually happens in the beginning.

Let $\lambda(\nu)$ be the maximal hyperbolic distance between the points
in the hyperbolic plane enclosed by an annulus of modulus $\nu$.
Note that $\lambda(\nu)=O(e^{-\nu})$ as $\nu\to\infty$.
Set $\lambda=\lambda(\mu)$.
 
Let $$h:(\Delta,V_i)\rightarrow (\tilde \Delta, \tilde V_i) \eqno (24)$$
be a $K$-qc pseudo-conjugacy between the corresponding configurations.
Our goal is to pull this map back to the next level.
The problem is that $h$ does not respect the positions of the critical
values. We assume first that we have a non-central return on this level,
that is $c_1\equiv g(c)\in V_k$ with $i\not=0$.

Let  $P_l$ be the pull-backs of $V_0$ by the univalent branches of 
iterates $g$
intersecting the critical set. We can  pull $h$ back by these
branches to obtain   a $K$-qc pseudo-conjugacy 
$$h_1 :  (\Delta, \cup P_l)\rightarrow (\tilde\Delta, \cup \tilde P_l).$$
This {\sl localizes} the positions of the critical values in the sense
that the hyperbolic distance between $h_1(c_1)$ and $\tilde c_1$ in $V_k$
is $O(\lambda)$. Indeed, they belong to the same puzzle-piece $\tilde P_l$
whose hyperbolic diameter in $V_k$ is at most  $\lambda$.

 Hence we can find a diffeomorphism
$\psi:\tilde \Delta\rightarrow\tilde\Delta$ which is id outside $\tilde V_k$,
moves $h_1 c$ to $\tilde c$, and has  a qc dilatation $1+O(\lambda)$.
Then the $K(1+O(\lambda))$-qc map
$$h_2=\psi\circ h_1: (\Delta, V_i)\rightarrow (\Delta, V_i)$$ 
respects the same configurations as $h$,
 and also carries $c$ to $\tilde c$.

Now we can pull $h_1$ back to 
$$H: (\Delta',U_i')\rightarrow (\tilde\Delta',\tilde U_i'), \eqno (25)$$
where $U_i'$ are $g$-pull-backs of $V_i$ to $\Delta'$. 
 However $U_i'$ are not the same as
$V_j'$, so we have to do more. 
What we need is to localize the
positions of the critical values  $a=g'c$ and $\tilde a$ of the next 
renormalizations. The argument depends on where they are. Let 
$a_1=g(a)\in V_j$. 

\smallskip {\sl Case (i).}  Assume $V_j$ is non-critical and different 
   from $V_k$. Then we can simultaneously move  of $c_1$
   and $a_1$ to the right positions, and then pull the map back to $\Delta'$.

\smallskip{\sl Case (ii).} Assume that $V_j=V_k$.
  If the hyperbolic distance between $a_1$ and $c_1$ in $V_j$ is greater
  than $\lambda(\mu/2)$ then the hyperbolic distance between the
  corresponding tilde-points is greater than 
  $\lambda(K\mu/2)$. Then we can simultaneously move these points to
 the right 
  positions by a qc map $\psi$   with dilatation $1+O(\lambda(K\mu/2))$.

  Otherwise let us first move  $c_1$ to the right position, and  pull the map
 back to $H$ as in (25). Then $a$ and $c$ stay in $V_0$ on
 hyperbolic distance $O(\lambda(K\mu/4))$, and the corresponding tilde-points
 stay on distance $O(\lambda(K^2\mu/4))$.  Hence $H(a)$ and $\tilde a$
  stay in $\tilde V_0$ on hyperbolic distance  $\delta=O(\lambda(K^2\mu/4))$.
  So we can move these points to the right positions by a $K(1+\delta)$-qc map
  respecting the boundary marking of $V_0$
  (though not respecting the critical points any more). 

\smallskip{\sl Case (iii).} Let us finally assume that $V_j=V_0$ is critical.
  Then $a$ belongs to a pre-critical puzzle-piece $V'_s\subset V_0$.
  Since mod($V_0\backslash V'_s)\leq \mu/2$, the map
   $H$ constructed above (see (25)) almost respects the positions of
  $a$-points in $V_0$.
  So we can make it respect these points keeping $\partial V_0$
   untouched.\smallskip 

  After all, we have constructed a $(1+O(\lambda^d))$-qc map 
$$h_2:(\Delta,V_i,a)\rightarrow (\tilde \Delta, \tilde V_i, \tilde a).$$

Let us now start over again, and  pull this map back to get a qc map
$$h_3 :  (\Delta, \cup P_l)\rightarrow (\tilde\Delta, \cup \tilde P_l).$$
 Unlike $h_1$ above this map also respects 
the forward orbits of the critical points 
(that is, the appropriate pre-images of $a$) until their
returns to the central puzzle-pieces. Hence
we can pull $h_3$ back to the critical puzzle-piece $\Delta'$.
Since $V_j'$ are
the pull-backs of $P_l$, this map respects the boundary marking of $V_j'$,
and we are done (in the non-central case).

\medskip{\bf Cascades of central returns.} Let
$V^m\supset V^{m+1}\supset...\supset V^{m+N}$
be a cascade of central returns, that is the critical value
$g_{m+1} c$ belongs to $V_k,\; k=m+1,...,m+N-1$, but escapes $V^{m+N}$.
We assume that the following configurations are qc pseudo-conjugate:
$$ h: (V^m, V^{m+1})\rightarrow (\tilde V^m, \tilde V^{m+1}).$$ 
Set $g=g_{m+1}, \; \mu={\rm mod}(V^m\backslash V^{m+1})$.

Let us take non-central puzzle-pieces 
$V^{m+1}_i\subset A^{m+1}=V^m\backslash V^{m+1}$
 and pull them back to the annuli $A^{m+2},...,A^{m+N}.$ We obtain a
Markov family of puzzle-pieces $W^{m+k}_i$. Let us induce on this
Markov scheme  {\sl the first landing map} 
$$\phi: \cup P^{m+k}_l\rightarrow V^{m+N}. $$ 
Then $${\rm mod}(W^{m+k}_i\backslash P^{m+k}_l )\leq\mu, $$
so that the dynamically defined  points are well localized by this partition.

Now we can proceed along the lines of the Main Construction
just using  the following substitution: $W^{m+k}_i$ play the role of
$V_i$, $V^{m+N}$ plays the role of $\Delta'\equiv V'$.
So we pull $h$ back to 
$$h_1: (\Delta, \cup P_l, V^{m+N})\rightarrow
  (\tilde\Delta, \tilde\cup P_l, \tilde V^{m+N}), $$
correct this map to make it  respect the $g$-critical values and then
pull it back to $V^{m+N}$ as in (25):
$$H: (V^{m+N}, U_i)\rightarrow (\tilde V^{m+N}, \tilde U_i), $$
where $U_i$ are the pull-backs of $W^{m+N}_i$ and $V^{m+N}$.
Take now  the first return $b$ of the
critical point back to $V^{m+N}$,  and look at Cases (i), (ii), (iii)
of the Main Step.
The first two cases go in the same way as above. However, the last case is
different since the pre-critical puzzle-pieces $U_1$ and $U_2$
are not necessarily well inside of $V^{m+N}$. 

To take care of this problem
let us first consider the first landing map 
$$\psi: \cup Y_j\rightarrow V^{m+N+1}$$ 
from $U_1\cup U_2$ to $V^{m+N+1}$,
and pull $H$ back to the domain of $\psi$. Since the components $Y_j$
 of this domain are
well inside $U_s$ (namely
 mod$(U_s\backslash Y_j)\geq {\rm mod}(A^{m+N+1})\geq\mu$), this gives us
an appropriate localization of the $b$-points.

\medskip\noindent {\bf Initial Construction.} In the beginning we have
a cascade $\Omega^0\supset...\supset \Omega^N$ of central returns with degenerate
annuli. So we may not directly apply the above argument.
Below we use the notations of  Lemma 0. We start with a qc pseudo-conjugacy
respecting the dynamics on the external rays through $\alpha$, $\alpha'$,
$\gamma$, $\gamma'$ and the equipotentials of $\partial \Omega^1$. 

\smallskip
{\sl Step 1.} Let us now construct a qc map $Q\rightarrow \tilde Q$.
Let $a=f^s c$ be the last point of orb$(c)$ landing at a $W_j$
before the return back to $V^0$.
 If the points $a=g^N c$ and $\tilde  a$
are well inside $W_j$ and $\tilde W_j$ correspondingly, then we can take
a qc map $(W_j,a)\rightarrow (\tilde W_j,\tilde a)$ and pull it back to
$Q$-pieces.   

Otherwise let us cut $\Omega^N$ by the external rays landing at 
$\gamma'$, and take the component $E$ attached to the fixed point $\alpha$.
Then the branch of $g^{-1}$ fixing $\alpha$ univalently maps $E$ into
itself. So $F=E\backslash g^{-1} E$ is a combinatorially well-defined 
fundamental domain for $g$ near the fixed point $\alpha$. Hence 
if $w\equiv f^i(a)\in E$ (combinatorially close to $\alpha$) then
there is the first moment $l\geq 0$ depending only on 
combinatorics such that $g^l w\in F$.

Let us also consider the fundamental domain 
$$F_*=F\cup g^{-1} F\cup g^{-2}F$$ for the 
third iterate of $g$.
If the external class of $f$ belongs to the given set of truncated limbs of
order two then the configuration  $(F_*,g^la)$ has a bounded geometry.
 Hence we can 
start with a qc map respecting these configurations and the dynamical pairing
on the $\partial F_*$.

Let us now pull this qc map back to $g^{-3}E,\; g^{-6}E,...$
If $l-1=3m$ then the $m$-fold pull-back will carry the point $w$ 
to $\tilde w$. Then we can pull this map back to $Q$-pieces by the
appropriate iterate of $f$.
If $l-1$ is not a multiple of 3 then we can replace
it by $l-1-n$ which is a multiple of 3 (where $n=1$ or 2) and correspondingly
replace $F$ by $g^{-n} F$.

\smallskip{\sl Step 2.} Let us now take a point $z\in \Omega^0$ 
and push it forward
by iterates of $g$ until it lands either at 
$\cup W_j\cup (\Omega^0\backslash \Omega^1)$
or at $Q$. If it happens, 
then we can pull the pseudo-conjugacy to an appropriate piece
containing $z$. The qc maps constructed in such a way agree on the common 
boundaries of the puzzle-pieces. The set of points where this map is not
defined is an expanding Cantor repellor. Hence it is qc removable, and
the  map automatically allows a qc continuation across it.  This provides us
with a qc map 
$$(\Omega^0, V^0, Q)\rightarrow (\tilde \Omega^0, \tilde V^0, \tilde Q)$$.

\smallskip{\sl Step 3.} Let us now localize the first return $b$
 of the critical  point back to $V^0$. To this end 
let us push $b$ forward 
 until the first moment $t$ it returns back to $Q$. Let $u=f^t b$.
The procedure depends on whether $u$ and $c$ stay on
  bounded hyperbolic distance in $Q$ (in terms of a given quantifier $R$)
  or not (compare Case (ii) of the Main 
  Construction). In the former case the position of $u$ is already 
   well localized inside of $Q$.
   In the latter case  we can on Step 1 simultaneously localize
   positions of $a=g^N c$ and $g^N u$ in $W_j$ and pull them back to $Q$.  
 Hence  we  can change the qc pseudo-conjugacy inside of 
  $Q$ so that it respects $u$-points. Pulling this back as on Step 2, 
  we construct a qc pseudo-conjugacy respecting $b$-points.

\smallskip {\sl Step 4.} Let us consider the full first  return map
   $G: \cup Z_j\rightarrow V^0$. Its domain covers the whole piece $V^0$
  except for a removable Cantor set $K\subset V^0$. 
  We can now construct a qs pseudo-conjugacy
  $$(V^0, Z_j)\rightarrow (\tilde V^0, \tilde Z_j) $$
  by the simple pull back and removing $K$. Since  the puzzle-pieces
  $V^1_i$ are among $Z_i$, we are done.

By Lemma 2, the principle modulus is definite on this level. 
So we can proceed further by applying the Main Step. 

\medskip\noindent{\bf Qc conjugacy on the critical sets.}
Let us show now that there is a qc map conjugating $f$ and $\tilde f$
on their critical sets.  Let $t=(m,n)\in T$ runs over the indices of the
principle nest of puzzle-pieces. Clearly the lexicographic order on $T$
corresponds to inclusion of the puzzle-pieces.
 Let $Q^t_0\equiv V^t$, and $Q^t_l$
be all pull-backs of  $Q^t_0$ corresponding to the first  landing
of the orbits of $z$, $z\in \omega(c),$ at $Q^t_0$. 
Then

$$ \omega(c)=\bigcap_t \bigcup_l  Q^t_l.$$
  Let us consider
the multiply connected domains 
$$P^t_l=Q^t_l\backslash\bigcup_k Q^{\tau}_k,$$
where $\tau\in T$ immediately follows $t$ in the lexicographic order.
The boundaries of $P^t_l$ are naturally marked.

By the {\sl Teicm\"{u}ller distance} between two marked domains
(of the same qc type) we mean the $\log K$ where $K$ is the qc
dilatation of the best qc homeomorphism between the domains respecting
the marking. 

\proclaim Lemma 11. The domains $P^t_l$ and $\tilde P^t_l$ stay on bounded
  Teichm\"{u}ller distance.

\noindent{\bf Proof.} We have proved that the pairs $(V^t, \cup_j V^{\tau}_j)$
and $(\tilde V^t, \cup_j \tilde V^{\tau}_j)$
stay a bounded Teichm\"{u}ller distance.  Pulling the corresponding qc
equivalence back by the univalent branches of $g_{\tau}$ we obtain that
$P^{\tau}_0$ and $\tilde P^{\tau}_0$ also stay a bounded Teichm\"{u}ller
distance. Pulling this back by the univalent branches of $f$ we obtain the
claim for all $l$.  \QED 

Gluing now together the multiply connected domains under consideration,
we construct a homeomorphism $h: V^0\rightarrow \tilde V^0$ which is qc on 
$V^0\backslash\omega(c)$ and conjugates $f$ and $\tilde f$ on their critical sets.
Since the critical sets are removable, we are done. 

\medskip\noindent{\bf Homotopy.} Let 
$\psi_0: (U,U')\rightarrow (\tilde U,\tilde U')$ be a homeomorphism
conjugating $f$ and $\tilde f$. We will show now that the qc map $h$
conjugating $f$ and $\tilde f$ on their critical sets can be constructed
in such a way that it is homotopic to $\psi_0$ rel the critical sets.

As in the proof of Corollary 2,
let $U_m\equiv V^{m,t(m)}$, $U_m'\equiv V^{m,t(m)+1}$, 
$U_m''\equiv V^{m,t(m)+2}$,
$G_m : U_m'\rightarrow U_m$ be the corresponding quadratic-like
renormalization of $f$, and let ${\bf Q}_m$ be defined as in (23).
These sets nest down to $\omega(c)$.

 A selection of the straightenings of
the quadratic-like  maps $G_m$ and $\tilde G_m$ provides us with
a choice of conjugacies 
$$\psi_m: (U_m',U_m'')\rightarrow (\tilde U_m', \tilde U_m'').$$
Let us continue $\psi_m$ to the annuli $U_m\backslash U_m'$ in such a way
that $\psi_m\simeq \psi_{m-1}$ (are homotopic) in the annulus 
$U_m\backslash J(G_m)$ rel the boundary.
Then let us spread $\psi_m$ around to the whole set ${\bf Q}_m$.
Outside ${\bf Q}_m$ set $\psi_m=\psi_{m-1}$. Clearly $\psi_m\simeq \psi_{m-1}$
mod $J_m$ where $J_m$ is the orbit of $J(G_m)$.

Let us define a homeomorphism $\psi : U\rightarrow \tilde U$ as the pointwise
$\lim\psi_m$. This homeomorphism is homotopic $\psi_0$ rel the critical sets.
Let us now construct a qc map $h$ homotopic to $\psi$ rel the critical sets.
First of all, the above selection of the  straightenings should be uniformly qc
which is possible because of the a priori bounds (Theorem B). 
 Then let us assume by
induction that we have already constructed a map $h_{m-1}\simeq \psi$ mod
${\bf Q}_m$ which is qc outside ${\bf Q}_m$. 

Let us  cut $U_m$ by the external rays through the points $\alpha, \alpha',
\gamma, \gamma'$ into puzzle-pieces $S_i$, and go through the above pull-back
construction. In the beginning we change $\psi_m$ on the $S_i$ to make it qc.
As the pieces $S_i$ are simply-connected, 
this change can be done via homotopy rel the boundary.
Then this homotopy can be pulled back to the deeper puzzle-pieces
according to the Starting Construction. This provides us with a 
homotopy rel the boundary
$$(V^0,\cup V^1_i)\rightarrow (\tilde V^0, \tilde V^1_0).$$
Then this homotopy can be pulled back through the cascade of Main Steps,
and spread around to the whole critical set (as in the previous subsection).
This gives us a qc map 
$$h_m: U'\backslash {\bf Q}_{m+1}\rightarrow \tilde U'\backslash \tilde {\bf Q}_{m+1}$$
homotopic to $\psi_m$ rel the boundary.

We should now continue this map to the annulus $R_m=U_m\backslash U_m'$.
To this end observe that $\psi_{m-1}$ has a bounded twist in this annulus
since it can be deformed rel the boundary to a qc map 
(by the above  pull-back argument).
Hence $\psi_{m-1}$ has a bounded twist in the annulus $U_m\backslash J(G_m)$
 as well, since this homotopy can be pulled back to this annulus (and 
by a hyperbolic argument will
automatically  be trivial on the Julia set).
Consequently the continuation of $\psi_m$ (and hence $h_m$) to $R_m$ (such that
$\psi_m\simeq \psi_{m-1}$ in $U_m\backslash J(G_m)$ mod the boundary)
 has a bounded twist as well. Hence this continuation can be realized
quasi-conformally.

Finally we can spread the homotopy from $U_m$ around the ${\bf Q}_m$.

\medskip\noindent
{\bf Sullivan's pull-back argument.} Remember that 
$\psi_0: (U,U')\rightarrow (\tilde U, \tilde U')$ is a 
conjugacy between $f$ and $\tilde f$, and $h$ is a $K$-qc map homotopic to
$\psi_0$ rel the critical sets. Sullivan's Pull-back argument allows us
to reconstruct $h$ into a qc conjugacy.

Let $U^n$ be the preimages of $U$ under the iterates of $f$. 
We can always assume 
that $h|U^n=\psi$.  Since $h(c_1)=\tilde c_1$, we can lift $h$ to a $K$-qc map
$h_1: U^1=\tilde U^1$ homotopic to $\psi$ rel the critical set and
$\partial U^1$. Hence $h_1=h$ on these sets, and we can continue $h_1$ to
$U\backslash U^1$ as $h$. This map conjugates $f$ and $\tilde f$ on 
the critical sets and also on $U^1\backslash U^2$.

Let us now replace $h$ with $h_1$ and repeat the procedure. In such a way
we construct a sequence of $K$-qc maps $h_n$ conjugating $f$ and $\tilde f$
on the critical sets and on $U^1\backslash U^{n+1}$. 
Passing to a limit we obtain a desired qc conjugacy.

\bigskip\centerline{\bf \S 5. Real case.}\medskip

In this section we will prove the following dichotomy:
 real maps of Epstein class (see below) either have a big complex space 
on the next quadratic-like level, or
essentially  bounded  real geometry (``essentially" loosely means
``up to saddle-node cascades"). The main ingredient is to create 
a generalized polynomial-like map with a definite modulus on an
essentially bounded level. By Theorem I  this implies big space,
provided the type is sufficiently high. From this dichotomy
we derive the real rigidity theorem.

%
 
\medskip\noindent{\bf Preliminaries.}
Let $\phi(z)=(z-c)^2$ denote the purely quadratic map. Let $I'\subset I$
be two nested intervals. A map 
$f: I'\rightarrow I$ is called {\sl quasi-quadratic} if it is $S$-unimodal
and has quadratic-like critical point $c$.  
 
Let us also consider a more general class ${\cal A}$ of maps  
$g: \cup T_i\rightarrow T$ defined
on a finite union of disjoint intervals $T_i$ compactly contained in an
interval $T$. Moreover, $g|T_i$ is a diffeomorphism onto $T$ for $i\not=0$,
while $g|T_0$ is unimodal with $g(\partial T_0)\subset\partial T$.
We also assume that the critical point $c\in T_0$ is quadratic-like, 
and that $S g<0$.  Maps of class ${\cal A}$ are real counterparts of 
generalized polynomial-like maps.  

Let $g\in {\cal A}$, and
$g|T_0=h\circ\phi$ where 
 and $h$ is a diffeomorphism of an appropriate interval 
$K\supset \phi(T_0)$ onto $T$. 
This map  belongs to the so-called {\sl Epstein class}
 ${\cal E}$ (see [S] and [L2])
if  the  inverse branches
$f^{-1}: T\rightarrow T_i$ for $i\not=0$  and $h^{-1}: T\rightarrow K$
allow analytic extension
 to the  slit complex plane 
${\bf C}\backslash T$. 


Let $I^0=[\alpha, \alpha']$ be the interval
between the dividing fixed point $\alpha$ and the symmetric one.
Let ${\cal M}$ denote the full  Markov
family of pull-backs of the interval $I^0$. Given a critical interval
$J\in {\cal M}$ (that is, $J\ni c$),
 we can define a (generalized)
renormalization $R_J f$ on $J$ 
as the first retun map to $J$ restricted to the components of its domain
 meeting the critical set. If $f$ admits a unimodal renormalization, then
 there are only finitely many
such components, so that we have a map of class ${\cal A}$. Moreover,
if $f$ is a map of Epstein class or a polynomial-like map, the
renormalization $R_J f$ inherits the corresponding property.


Let $I^0\supset I^1\supset\ldots\supset I^{t+1} $ be the real principal nest
of intervals until the next quadratic-like level (that is,  $I^{n+1}$ is the
pull-back of $I^n$ corresponding to the first return of the critical point.
Let us use the same notation
 $g_n: \cup I^n_j\rightarrow I^{n-1}$ for the real
generalized renormalizations on the intervals $I^n$.

Our first goal is to fill-in the gap in between the notions of bounded 
combinatorial type
in the sense of period and in the sense of the number of central cascades.
To this end we need to analyse in more detail cascades of central returns.

The return on level $n-1$ is called {\sl high} or {\sl low} if
$g_n I^n\supset I^{n}$  or $g_n I^n\cap I^n=\emptyset$ correspondingly.
Let us classify the cascades
$$I^m\supset...\supset I^{m+N},\;\;  g_{m+1}c\in I^{m+N-1}\backslash I^{m+N}
\eqno (26)$$
 of central returns
as {\sl Ulam-Neumann} or {\sl saddle-node} according as the return on the
level $m+N-1$ is high or low. There is a fundamental difference between these
two types of cascades. Let us call the levels $m+1,...,m+N-1$ 
of  a saddle-node cascade {\sl neglectable}, and all other levels 
{\sl essential}. Let $m=e(l)$ counts the essential levels.

Let $K^{m+i}_j\subset I^{m+i-1}\backslash I^{m+i}$ denote
the pull-back of $I^{m+1}_j$ under $g_{m+1}^{\circ(i-1)}$, $i=1,...,N$, $j\not=0$.
Clearly  $K^{m+i+1}_j$ are mapped by $g_{m+1}$ onto $K^{m+i}$, $i=1,..., N-1$,
while $K^{m+1}_j\equiv I^{m+1}_j$ are mapped onto the whole $I^m$.
So we have a Markov scheme associated with any central cascade.


Take now a point $x\in \omega(c)\cap (I^m\backslash I^{m+1})$ on an
essential level $m=e(l)$. 
Let us push it forward 
by iterates of $g=g_{m+1}$ through the above Markov scheme until it lands
at the next essential level $I^{m+N}$, $m+N=e(l+1)$.
 Let $o(x)$ (``the order of $x$")
 denote the number of times it passes through
 $I^m\backslash I^{m+1}$ before landing at $I^{m+N}$
(e.g., $o(x)=1$ if $g x\in I^{m+N}$). Let $g x\in I^{m+i}$.
Then set $d(x)\min \{i, N-i\}$ (``the depth of the first iterate").

Let us now introduce the scaling factors
$$\lambda_n\equiv \lambda_n(f)={|I^n|\over |I^{n-1}|}.$$
According to [L2],
these scaling factors exponentially decay with the number
of central cascades.
Moreover, this rate is uniform when the scaling factors become small enough.

Let us call the geometry of $f$  {\sl essentially bounded} 
(until the next quadratic-like level)
if the scaling factors $\lambda_n=|I^n|/|I^{n+1}|$ stay away from 0,
while the configurations $(I^{n-1}\backslash I^n, I^n_k)$ have
bounded geometry (that is, all intervals $I^n_j,\; j\not=0$, and all components of
$I^{n-1}\backslash \cup I^n_k$  (``gaps")
are commensurable). Remark that we allow the scaling factors $\lambda_n$
to be close to 1.

\medskip
\noindent{\bf Complex bounds.}
 Sullivan's Sector Lemma provides us with complex
bounds in the case when $f$ is infinitely q-renormalizable of bounded type.
In the  non-q-renormalizable case 
 the complex bounds were obtained in [LM] and [L2].
We will complement these results with the following theorem.
 
Let us pick a class ${\cal U}_{\tau, \bar\mu}$ 
  of real quadratic-like maps $f$
  of the same q-renormalizable type $\tau$, and such that
   mod$(f)\geq \bar\mu$, where $\bar\mu>0$
 is an a priori chosen small quantifier.

\proclaim Theorem D.   One of the following two possibilities occurs
  for all $f\in {\cal U}_{\tau,\mu}$ simultaneously:
  either mod$(R f)\geq\bar\mu>0$, or
  the real geometry  of $f$ is essentially bounded (until the next
quadratic-like level).
 

In the following two lemmas we analyse the geometry of long central cascades.
Let us call a quasi-quadratic map 
saddle-node or Ulam-Neumann if it is topologically
conjugate to $z\mapsto z^2+1/4$ or $z\mapsto z^2-2$ correspondingly. 

 \proclaim Lemma 12. Let us consider an Ulam-Neumann cascade as (26)
 with commensurable $I^m$ and $I^{m+1}$. Then there is a bounded 
  $l$ such that the generalized  renormalization $g_{m+l}$ allows
  a polynomial-like extension to the complex plane with a definite 
  modulus. Moreover, the principle modulus $\mu_{m+N+1}$
  is big, provided the cascade is long.

\noindent{\bf Proof.}  It is easy to see by compactness argument
that if the cascade is long
enough then  the map 
$g_{m+1}: I^{m+1}\rightarrow I^m$ (with the domain rescaled to the
unit size) is $C^1$-close to an Ulam-Neumann map. It follows that
$I^{m+N}$ occupies a definite part of $I^m$, and, moreover,
$|I^{m+k}\backslash I^{m+N}|$  decrease with $k$
at a uniformly exponential rate. 
Hence there is a bounded $l$ such that $I^{m+l}\backslash I^{m+N}$
is $\epsilon$-tiny as compared with $I^{m+N}$. 

Take now the Euclidian disk $D=D(I^{m+l})$ and pull it back by the
inverse branches of $g_{m+l+1}$ (as in the previous lemma). 
As $g_{m+l+l}=h\circ \phi$ where $h$ is a diffeomorphism with a
bounded distortion, the central pull-back will be an ellipse based
upon the interval $I^{m+l+1}$
whose imaginary axis is $O(\sqrt{\epsilon}|I^{m+l+1}|)$. It follows that
this ellipse is well inside of $D$.  

 
The last statement follows from Lemma 11.
\QED





\proclaim Lemma 13. All saddle-node patterns (26) of the same length
  with commensurable 
 $I^m$ and $I^{m+1}$ are qs equivalent.

\noindent{\bf Proof.} Let 
 $g: I'\rightarrow [0,1]$ be a unimodal map
 of Epstein class (and perhaps escaping critical point): $g\in {\cal E}_u$.
 By definition, $g=h\circ\phi$ with a diffeomorphism $h$ whose 
inverse allows the analytic extension to ${\bf C}\backslash [0,1]$.
Let us supply this space with   with the Montel topology on $h^{-1}$.

The set of $g \in {\cal E}_u$ with bounded geometry on the real line
is compact. Hence given a long saddle-node cascade (26), the map $G$
obtained from 
 $g_{m+1}: I^{m+1}\rightarrow I^m$ by rescaling $I^m$
 to the unit size must be close
 to a saddle-node quadratic-like map.
Hence we can reduce $G$
to a form $z\mapsto z+\epsilon+\psi(z)>0$ where $\psi(z)$ is uniformly
comparable with $z^2$, and (as we will see in a moment)
 $\epsilon$ is determined (up to a bounded error) by the length of the
cascade.

Take a big $a>0$. When $|z|<a\sqrt\epsilon$, the step $G(z)-z$ is 
of order $\epsilon$. Otherwise $\psi(z)$ dominates over $\epsilon$,
and in the chart $\zeta=1/z$ the step is of order 1. It follows that
the qs class of the cascade is determined by $\epsilon$, which in turn 
is related to the  length of the cascade  by
$N\asymp 1/\sqrt{\epsilon}$.     \QED

Given two intervals $L\subset S$,
let $P(L|S)$ denote the Poincar\'{e}  length of $L$ in $S$. 
Given an interval $I$,
let $D(I)$ denote the Euclidian disk based upon $I$ as a diameter.

\proclaim Lemma 14.
Assume that $\lambda_n<\epsilon$ with a sufficiently small $\epsilon>0$. 
Then there is an interval $T\in {\cal M}$ containing $I^{n+2}$ 
such that the renormalization $R_T f$
allows a (generalized) polynomial-like extension
to the complex plane with the
principle modulus $\mu\to\infty$ as $\epsilon\to 0$.

\noindent{\bf Proof.} First of all we can assume that all intervals $I^n_j$
are well inside $I^{n-1}$ (otherwise pass to the next level). 
The following  construction of a polynomial-like map
is combinatorially the same  
as in [L], Lemma 5.3. For the reader's convenience we briefly repeat it.

Let $g\equiv g_n: \cup I^n_j\rightarrow I^{n-1}$.  Let us inductively define
the {\sl cut-off orbit} of $I^n_0$ as
$$g^l_{cut} (I^n_0)= g\,(g^{l-1}_{cut}(I^n_0)\cap I^n_j),$$
provided $I_j\ni g^{l-1}c, \; j\not=0$. We stop at the first  moment when
$g^l_{cut} I^n_0\cap I^n_0\not=\emptyset$. Let us define  
$T\equiv T_0\ni c$  as the 
pull-back of $I^{n-1}$ by $g^l$, and set $G|T_0=g^l$. Clearly
$I^n\supset T_0\supset I^{n+1}$, and  $G(\partial T_0)\supset\partial I^{n-1}.$

Let now $z\in (\omega(c)\cap I^{n-1})\backslash T_0$.  
If $z\in I^{n-1}\backslash I^n$ then let $T(z)$ be the interval 
$I^n_j\equiv I^n(z)$ containing $z$, and $G|T(z)=g$.
If $z\in I^n\backslash T_0$ then let us push $z$ forward by iterates of $g$
until it is separated from the corresponding iterates of $c$ by the intervals
$I^n_j$. Let it happen at moment $s$, and $g^s c\in I^n_k$. It follows from
the choice of $l$ that $s\leq l$ and $k\not=0$. Let us now define
$T(z)$ as the pull-back of $I^n_k$ by $g^s$, and $G|T(z)=g^{s+1}$.
It is easy to see that $G: T(z)\rightarrow I^{n-1}$ is a diffeomorphism.

So we have constructed a map $G: \cup T_i\rightarrow I^{n-1}$ 
of class ${\cal A}$.
Let us now take the Euclidian disk $D=D(I^{n-1})$ and pull it back by
the inverse branches of $G$. This provides us with a set of domains $D_i$
based upon the intervals $T_i$. Moreover, by a little hyperbolic
argument (see e.g., Lemma 8.1 of [LM]) $D_i\subset D(T_i)$.

Let us now estimate the shape of $D_0$. To this end let us consider
the  following decomposition:
$$G|T_0=(g|I^n_i)\circ (h|K)\circ\phi|T_0. $$
Here  $g^{l-1}c\in I^n_i$, and $g^{l-1}|T_0=h\circ\phi$,
where $h$ is a diffeomorphism of an appropriate interval $K$ onto  $I^n_i$
with a Koebe space spreading over $I^{n-1}$. 
As $I^n_i$ is well inside of $I^{n-1}$, $h|K$ has a bounded distortion.
Moreover, $g|I^n_i$ is quasi-symmetric (as a composition of the quadratic map
and a diffeomorphism of bounded distortion). Hence
$G|T_0=(H|K)\circ\phi|T_0$
with a quasi-symmetric diffeomorphism $H$. Furthermore, as
 $G(T_0)\cap I^n_0\not=\emptyset$, 
$$|G(T_0)|\geq {1-\epsilon\over 2} |I^{n-1}|.$$
Pulling this back by the qs map $H$, we conclude that
$$|\phi(T_0)|\geq\delta|K|$$
with $\delta=\delta(\epsilon)$. Let $Q$ be the pull-back of $D$ by $H$.
Then $Q\subset D(K)$. Pulling this back by the quadratic map $\phi$,
we conclude that $D_0$ has a bounded shape. As it is based upon a
$\epsilon$-tiny interval $I^n_0$, it is well inside $D$. Moreover,
the annulus $D\backslash D_0$ is getting big as $\epsilon\to 0$. 

It follows that $R_T G$ satisfies the desired properties. 
Finally, it is easily seen from the construction that the first
 return map to $T$ under $f$ coincides
with the first return map under $f$, so that $R_T f=R_T G$. \QED 

Now we are ready to state the key lemma.
 
\proclaim Lemma 15.
There is an interval $T\in {\cal M}$  
such that the renormalization $R_T f$
allows a polynomial-like continuation to the complex plane with a
definite principle modulus $\mu$. Moreover, $T$ lies on an essentially bounded
level:   $T\supset I^{e(l)}$.

 
\noindent{\bf Proof.}
 Take a small $\epsilon>0$ and $\delta>0$,
 and select the first moment $l$ for which
$$\lambda_{l}>(1-\delta)\lambda_{l-1}. \eqno (27)$$
For such a level [L2,\S 5] provides us with a polynomial-like map
$G: \cup D_i\rightarrow D(I^l)$ with a definite modulus and such that
 the number of central cascades preceding $T_0=D_0\cap{\bf R}$ 
is bounded. Moreover, only the last of these cascades may be of 
Ulam-Neumann type. If this cascade is of bounded length then $T_0$
lies on an essentially bounded level. Otherwise
 Lemma 12 provides us with a 
desired polynomial-like map.

On the other hand, if (27) fails to happen on the first
$s=\log\epsilon/log(1-\delta)+1$ levels then we come up with an
$\epsilon$-small scaling factor, and can  apply Lemma 14.  \QED

Given a q-renormalizable map $f$, let $\tau(f)$ denote the maximum of the 
type $\kappa(f)$, the lengths of the Ulam-Neumann cascades, 
the orders $o(x)$ and the depths $d(x)$ for all $x\in \omega(c)$.

\proclaim Lemma 16. Take  a  $\bar\mu>0$.
  Let $f$ be a q-renormalizable unimodal map of Epstein class
  of a bounded distortion $D$. If $\tau (f)$ is sufficiently high
  (depending on $D$ and $\bar\mu$ only), then 
the renormalization $Rf$ is polynomial-like with mod$(Rf)>\bar\mu$.  




\noindent{\bf Proof.} 
 Assume that (i)  occurs.
 Then  by Lemma 15 on an essentially bounded level $T$ we can
  create a generalized 
  polynomial-like map  $R_T f$ with a definite modulus 
   ($>\bar \nu>0$). Then by Theorem~A the moduli of further renormalizations
   of $R_T f$ will grow at a linear rate with the number of central cascades.
    Hence the quadratic-like renormalization $Rf$ will have 
    a $\bar\mu$-big modulus, provided there are sufficiently many central
  cascades.
 
If (ii) occurs then by Lemma 12 in the end of the Ulam-Neumann cascade 
we observe a generalized polynomial-like map with a big modulus.
Then by Corollary 6 the modulus
of the quadratic-like renormalization $Rf$ will be big as well.


  Assume further that there is an $x\in \omega(c)\cap (I^m\backslash I^{m+1})$
of a high order $o(x)$,
where $I^m\supset...\supset I^{m+N}$ is a central cascade as (26)
(it may be $N=1$). Let us consider the above Markov scheme involving the
 intervals $K^{m+i}_j$. Let $J\ni x$ denote the pull-back of $I^{m+N}$
 corresponding to the first landing of the orb$(x)$ at $I^{m+N}$.

As the intervals $K^{m+1}_j$ are well inside of 
  $I^{m}\backslash I^{m+1}$, and orb$(x)$ passes many times 
  through these intervals before the first landing
  at $I^{m+N}$, the Poincar\'{e} length \break
 $P(J|(I^m\backslash I^{m+1}))$ is big.
  Pulling this interval back
  to the critical point  we will find
  a level with a small scaling factor.   Applying Lemma 14 
  we get the claim. 

  Let us finally assume that there is an 
 $x\in \omega(c)\cap (I^m\backslash I^{m+1})$ with high $d(x)$. 
Then $g_{m+1}x\in I^{m+i}\backslash I^{m+i+1}$ with $d(x)\leq i\leq N-d(x)$.
Then by Lemma 13  $I^{m+i}\backslash I^{m+i+1}$ is tiny in $I^m$. It follows
 that the interval $J\ni x$ introduced two paragraphs up is tiny in
  $I^m\backslash I^{m+1}$. Now we can complete the argument as above. 
  \QED

\smallskip\noindent{\bf Remark.} 
Now a little extra work shows that if 
$\tau(R^mf)$ is sufficiently high on all levels, then MLC
holds at $c\in {\bf R}$.


\proclaim Lemma 17. If $\tau(f)$ is bounded, then
the geometry of $f$ is essentially bounded
 (until the next quadratic-like level).

\noindent{\bf Proof.} Assume that the geometry is bounded on level $n-1$, and
let us see what happens on the next level. Given an 
$x\in \omega(c)\cap (I^{n-1}\backslash I^n)$, let $J(x)$ denote the
pull-back of $I^{n}$ corresponding to the first landing of orb($x$) at
$I^n$. As the landing time under iterates of $g_n$ is bounded,
$J(x)$ is commensurable with $I^{n-1}$. 

To create the intervals $I^{n+1}_j$, we should pull all intervals $J(x)$
back by $g_n: I^n\rightarrow I^{n-1}$. As $g_n$ is a quasi-quadratic map,
all non-central intervals $I^{n+1}_j$ and the gaps in between 
 are commensurable with $I^n$. 

The only possible problem is that the central interval $I^{n+1}$ may be
tiny in $I^n$. This may happen only if the critical value $g_n c\in J(x)$
is very close to the $\partial J(x)$. Let $l$ be such that 
$f^l J(x)=I^n$. Since $f^l: J(x)\rightarrow I^n$ is qs, 
$g_{n+1}c=g_n^{\circ (l+1)}$ turns out to be very close to $\partial I^n$
(``very low return"). But $g_{n+1}c$ belongs
to some non-central interval $I^{n+1}_j$ whose Poincar\'{e} length in
$I^n$ is definite (as we have shown above). This is a contradiction.

  So when we pass from one level to the next, the  geometric bounds
change gradually (provided the conditions of Lemma 16 don't hold).
But the same is true when when we pass from level $m=e(l)$ to level
 $m+N=e(l+1)$
of a saddle-node cascade (26). Indeed, assume that the geometry on level
$I^m$ is bounded. Then the geometry of all configurations 
$(I^{m+i-1}\backslash I^{m+i}, K^{m+i}_j)$, $i=1,...,m+N$, are bounded as well.
Let us define the intervals $J(x)$,
 $x\in\omega(c)\cap (I^{m+N-1}\backslash I^{m+N})$, as 
the pull-backs of $I^{m+N}$ corresponding to the first landing of orb$(x)$
at $I^{m+N}$. Then it follows from boundedness of $o(x)$ and $d(x)$ that
the configurations of intervals $J(x)$ has a bounded geometry in $I^{m+N-1}$. 
Now we can pull these intervals  back to the next level $m+N$, and
argue that the geometry is still bounded in the same way as above.  \QED 

Now Theorem D follows from the last two lemmas.

\medskip\noindent{\bf Quasi-symmetric conjugacy.}
We will show below that any two real quadratic-like maps with the same
combinatorics are qs conjugate, which implies the real rigidity
conjecture (compare [Sw]). To construct the conjugacy,
we bounce in between Sullivan's argument
for bounded geometry case, and the pull-back
argument of \S 4.
 
Let us take two maps $f$ and $\tilde f$ of Epstein class
with a  bounded distortion on the real line.
Let us consider the alternatives of Theorem D. In the latter case
  the real geometry is essentially bounded before the next quadratic-like level.
This allows us to construct a qc pseudo-conjugacy between the configurations
of the Euclidian disks based upon the intervals $I^n_j$.
 The construction is the same as in the bounded geometry case
(see [MS], Ch. IV, Theorem 3.1), 
except that Lemma 12 takes care of
 long saddle-node cascades.
 
 
If the first alternative of Theorem D occurs, then by Lemma 15
on some essentially bounded level
we can create  polynomial-like maps with definite moduli.
 By Lemma 17  the geometry is essentially bounded until that level, and
we can apply the previous argument. On that level we can switch 
to the pull-back argument of \S 4. (To begin the
argument, use the initial construction of [L3], \S 5.) 
 
When we arrive at the next quadratic-like level, then we proceed as follows.
In the first case we have arrived with a qc pseudo-conjugacy between
configurations of Euclidian disks. Then just apply the previous construction
to $Rf$ (here we need real a priori bounds for infinitely q-renormalizable maps
[G], [BL], [S]).
 In the second case we have arrived with configurations of
topological disks. Then  interpolate the qc pseudo-conjugacy as in \S 4,
and conformally map the range of $Rf$ to a slit domain. This gives us a
map of Epstein class  with a definite distortion on the real line,
and we can repeat the construction.

\bigskip
\centerline{\bf References.}\smallskip

{\item [AB]} L. Ahlfors \& A. Beurling. Conformal invariants and 
  function-theoretic null-sets. Acta Math. {\bf 83} (1950), 101-129.\smallskip

{\item [BL]} A.Blokh \& M.Lyubich. Measure and dimension of solenoidal
  attractors of one dimensional dynamical systems. Comm. Math. Phys.
  {\bf 127} (1990), 573-583.\smallskip

{\item [BH]} B.Branner \& J.H.Hubbard.   The iteration of cubic polynomials,
  Part II.  Acta Math. {\bf 169} (1992), 229-325.\smallskip

{\item [D1]} A. Douady. Shirurgie sur les applications holomorphes.
  In ``Proc. Internat. Congress Math. Berkeley", {\bf 1} (1986), 724-738.
  \smallskip

{\item [D2]} A. Douady. Description of compact sets in ${\bf C}$.
  In:
    ``Topological Methods in Modern Mathematics, A Symposium in
     Honor of John Milnor's 60th Birthday", Publish or Perish, 1993.\smallskip

{\item [DH1]} A.Douady \& J.H.Hubbard. \'{E}tude dynamique des polyn\^{o}mes
   complexes. Publication Mathematiques d'Orsay, 84-02 and 85-04.\smallskip

{\item [DH2]} A.Douady \& J.H.Hubbard. On the dynamics of polynomial-like maps,
    Ann. Sc. \'{E}c. Norm. Sup. {\bf 18} (1985), 287-343.\smallskip

{\item [F]} P. Fatou. M\'{e}moires sur les \'{e}quations fonctionnelles.
  Bull. Soc. Math. France, {\bf 48}, 33-94. 

{\item [G]} J. Guckenheimer. Limit sets of $S$-unimodal maps with zero
  entropy. Comm. Math. Phys. {\bf 110}, 655-659.\smallskip

{\item [H]} J.H. Hubbard.  Local connectivity of Julia sets and bifurcation
     loci: three theorems of J.-C. Yoccoz. In:
    ``Topological Methods in Modern Mathematics, A Symposium in
     Honor of John Milnor's 60th Birthday", Publish or Perish, 1993.\smallskip

{\item [HJ]} J. Hu \& Y. Jiang. The Julia set of the Feigenbaum quadratic
  polynomial is locally connected. Preprint 1993.\smallskip

{\item [J]} Y. Jiang. Infinitely renormalizable quadratic Julia sets.
  Preprint, 1993.\smallskip

{\item [K]} J. Kahn. Holomorphic Removability of Julia Sets. Manuscript in
  preparation since 1992.\smallskip

{\item [LM]} M. Lyubich \& J.Milnor. The unimodal Fibonacci map. 
         Journal of AMS, {\bf 6} (1993), 425-457. \smallskip


{\item [L1]} M. Lyubich. On the Lebesgue measure of the Julia set of a 
     quadratic polynomial, Preprint IMS at Stony Brook, 1991/10. \smallskip

{\item [L2]} M. Lyubich. Combinatorics, geometry and attractors of
     quasi-quadratic maps.    
     Preprint IMS  at Stony Brook \#1992/18.\smallskip

{\item [L3]} M. Lyubich. Teichm\"{u}ller space of Fibonacci maps.
  Preprint, 1993.\smallskip

{\item [MSS]} R. Ma\~{n}\'{e}, P.Sad, D. Sullivan. On the dynamics of rational
  maps. Ann. scient. Ec. Norm. Sup., {\bf 16}  (1983), 193-217.\smallskip

{\item [M1]} J. Milnor. Self-similarity and hairiness in the Mandelbrot set,
          pp. 211-257 of ``Computers in geometry and topology",
          Lect. Notes in Pure Appl Math, {\bf 114}, Dekker 1989.\smallskip

{\item [M2]} J. Milnor. Local connectivity of Julia sets: expository
   lectures. Preprint IMS Stony Brook, \#1992/11.\smallskip

{\item [McM]} C. McMullen. Complex dynamics and renormalization.
  Preprint, 1993.\smallskip

{\item [MS]} W. de Melo \& S. van Strien. One dimensional dynamics.
   Springer-Verlag, 1993. \smallskip

{\item [Sw]} G. Swiatek. Hyperbolicity is dense in the real quadratic
  family.  Preprint IMS  at Stony Brook \#1992/10.\smallskip

{\item [S]} D.Sullivan. Bounds, quadratic differentials, and renormalization
conjectures.  AMS Centennial Publications.
{\bf 2}: Mathematics into Twenty-first Century (1992). \smallskip

\end